\documentstyle[12pt]{article}
\newcommand{\sect}[1]{\section{#1}\setcounter{equation}{0}}

\font\mbn=msbm10 scaled \magstep1
\font\mbs=msbm7 scaled \magstep1
\font\mbss=msbm5 scaled \magstep1
\newfam\mbff
\textfont\mbff=\mbn
\scriptfont\mbff=\mbs
\scriptscriptfont\mbff=\mbss
\def\mbf{\fam\mbff}
\def\Re{{\mbf R}}

\def\Z{{\mbf Z}}
\def\Co{{\mbf C}}

\newtheorem{Th}{Theorem}[section]
\newtheorem{Lm}[Th]{Lemma}

\newtheorem{D}[Th]{Definition}
\newtheorem{Proposition}[Th]{Proposition}
\newtheorem{R}[Th]{Remark}
\author{Alexander Brudnyi\thanks
{Research supported in part by NSERC.
\protect\newline 
}\\
}
\title{CLASSIFICATION THEOREM FOR A CLASS OF FLAT CONNECTIONS AND
REPRESENTATIONS OF K\"{A}HLER GROUPS}
\date{November, 1998}
\begin{document}
\maketitle
\begin{abstract}
{The paper presents a classification theorem for the class of
flat connections with triangular (0,1)-components on a topologically trivial
complex vector bundle over a compact K\"{a}hler manifold. As a consequence 
we obtain several results on the structure of K\"{a}hler groups, i.e.,  
the fundamental groups of compact K\"{a}hler manifolds. 
} 
\end{abstract}
\sect{\hspace*{-1em}. Introduction.}
{\bf 1.1.} 
Let $M$ be a compact K\"{a}hler manifold. For a matrix Lie group $G$
the representation variety
${\cal M}_{G}$ of the fundamental group $\pi_{1}(M)$ is
determined as $\penalty-10000$
$Hom(\pi_{1}(M),G)/G$.  Here
$G$ acts on the set $Hom(\pi_{1}(M),G)$ by pointwise
conjugation: 
$(gf)(s)=gf(s)g^{-1}$,\ $s\in\pi_{1}(M)$. 
A study of geometric properties of ${\cal M}_{G}$ is of 
interest because of the relation to the problem of classification of 
K\"{a}hler groups (the problem was posed
by J.-P. Serre in the fifties). 
For a simply connected nilpotent Lie group $G$ every element of 
${\cal M}_{G}$ is uniquely determined by a $d$-harmonic nilpotent 
matrix 1-form $\omega$ on $M$ satisfying $\omega\wedge\omega$ represents
0 in the corresponding de Rham cohomology group.
It follows, e.g., from the theorem of Deligne-Griffiths-Morgan-Sullivan 
on formality of a compact K\"{a}hler manifold
(see [DGMS]).
The main result of our paper gives, in particular, a similar
description for elements of ${\cal M}_{G}$ with a simply connected
solvable Lie group $G$.
Our arguments are straightforward and based on cohomology techniques only.
As a consequence of the main theorem we obtain several results on the 
structure of K\"{a}hler groups. 
We now proceed to formulation of the results.

It is well known that ${\cal M}_{GL_{n}(\Co)}$ is equivalently 
characterised as moduli spaces of flat bundles over $M$ with structure
group $GL_{n}(\Co)$.
In this paper we deal with a family of $C^{\infty}$-trivial complex flat vector
bundles over $M$.
Every bundle from this family is determined by a flat connection on 
the trivial bundle $M\times\Co^{n}$, i.e. 
by a matrix valued 1-form $\omega$ on $M$ satisfying 
\begin{equation}\label{form2}
d\omega-\omega\wedge\omega=0.
\end{equation}
Moreover, we assume that the (0,1)-component $\omega_{2}$ of $\omega$
is an upper triangular  matrix form. Denote this class of connections
by ${\cal A}_{n}^{t}$. 
\begin{R}\label{conexam}
{\rm E.g., connections from ${\cal A}_{n}^{t}$  determine
(by iterated path integration) all representations of $\pi_{1}(M)$ into 
simply connected real solvable Lie groups. (Here according to Lie's 
theorem we think of every such group as a subgroup of a complex 
Lie group of upper triangular matrices.)}
\end{R} 

Let $T_{n}(\Co)$ denote the complex Lie group of upper triangular
matrices in $GL_{n}(\Co)$. Then
the group $C^{\infty}(M,T_{n}(\Co))$ acts by $d$-gauge transforms on the 
set ${\cal A}_{n}^{t}$:
\begin{equation}\label{form3}
d_{g}(\alpha)=g^{-1}\alpha g-g^{-1}dg,\ \ \ \ 
(g\in C^{\infty}(M,T_{n}(\Co)),\ \alpha\in {\cal A}_{n}^{t}).
\end{equation}
Denote the corresponding quotient space by ${\cal B}_{n}^{t}$.
In this paper we study the structure of ${\cal B}_{n}^{t}$.
Also our result gives a characterization of the subset of
${\cal M}_{GL_{n}(\Co)}$ consisting of conjugate classes of
representations determined by elements of ${\cal A}_{n}^{t}$.

Let ${\cal U}_{\oplus}^{n}$ be a class of flat vector bundles over $M$ of 
complex rank $n$ whose elements are direct sums of
topologically trivial flat vector bundles of complex rank 1 with unitary
structure group. Note that every $E\in {\cal U}_{\oplus}^{n}$
should be constructed by a unitary diagonal
cocycle $\{c_{ij}\}_{i,j\in I}$ defined on an open covering 
$\{U_{i}\}_{i\in I}$. All definitions formulated below do not depend
on the choice of such cocycle.\\
A family $\{\eta_{i}\}_{i\in I}$ of matrix-valued $p$-forms 
satisfying
\begin{equation}\label{section}
\eta_{j}=c_{ij}^{-1}\eta_{i}c_{ij}\ \ \ on\ \ \ U_{i}\cap U_{j}
\end{equation}
is, by definition, a $p$-form with values in the bundle $End (E)$.
We say that such form is {\em nilpotent} if every $\eta_{i}$ takes its 
values in the Lie algebra of the Lie group of upper triangular unipotent 
matrices.
Since $End(E)\in {\cal U}_{\oplus}^{n^{2}}$,
there exists a natural flat Hermitian metric on $End(E)$. Then, as usual,
we construct by this metric a $d$-Laplacian on the
space of $End(E)$-valued forms. In what follows harmonic forms are
determined by this Laplacian.
Denote by ${\bf H}_{d}^{1}(End(E))$ the finite-dimensional complex vector
space of $End(E)$-valued harmonic 1-forms and by $H^{2}(End(E))$ the de Rham 
cohomology group of $End(E)$-valued $d$-closed 2-forms. Further, 
consider the set ${\bf H}_{0}^{t}(End(E))\subset {\bf H}_{d}^{1}(End(E))$ of
harmonic forms $\eta$ satisfying
$$
\begin{array}{l}
(i)\ \ \ (0,1)-{\rm component}\ \eta_{2}\ {\rm of}\ \eta \
{\rm is\ nilpotent};\\
(ii)\ \ \eta\wedge\eta \ {\rm represents}\ 0\ {\rm in}\ H^{2}(End(E)).
\end{array}
$$
Observe that
${\bf H}_{0}^{t}(End(E))$ is a complex affine subvariety of 
${\bf H}_{d}^{1}(End(E))$ defined by homogeneous quadratic equations.

Let $Aut_{f}^{t}(E)$ be the group of triangular flat authomorphisms of 
$E$. Elements of $Aut_{f}^{t}(E)$ are, by definition, locally constant
sections of $End(E)$ determined by (\ref{section}) with
$\eta_{i}\in T_{n}(\Co)\  (i\in I)$. Clearly, $Aut_{f}^{t}(E)$ is a
complex solvable Lie group. It acts by conjugation on the space of
$End(E)$-valued forms and commutes with the Laplacian. In particular,
it acts on ${\bf H}_{0}^{t}(End(E))$. Consider the quotient space
${\cal S}_{E}^{n}:={\bf H}_{0}^{t}(End(E))/Aut_{f}^{t}(E)$ and denote by
${\cal S}^{n}$ the disjoint union 
$\sqcup_{E\in {\cal U}_{\oplus}^{n}}
{\cal S}_{E}^{n}$. (Note that according to Green-Lazarsfeld theorem [GL]
if the dimension of the image of the Albanese mapping of $M$ 
$\geq 2$ the set ${\cal S}_{E}^{n}$ with the generic $E$ consists of a single 
point.) 
\begin{Th}\label{maintheo}
There is a one-to-one correspondence between the sets ${\cal B}_{n}^{t}$
and ${\cal S}^{n}$.
\end{Th}
Using Theorem \ref{maintheo} for the case of flat connections corresponding
to unipotent representations of $\pi_{1}(M)$ one can give
alternative proofs of some results, e.g. due to Campana
(for references see [ABCKT]), Gordon and Benson [GB]. 
In the following section we describe the above correspondence in more 
details.
\\
{\bf 1.2.} We now formulate several geometrical applications of Theorem 
\ref{maintheo}. They describe some properties of the set $S_{n}(M)$ of 
representations of $\pi_{1}(M)$ into $GL_{n}(\Co)$ generated by 
connections from ${\cal A}_{n}^{t}$. 

Let $T_{2}^{u}$ denote the Lie group of upper-triangular $(2\times 2)$-
matrices with unitary elements on the diagonal. Further, denote by
$S_{2}^{u}(M)$ a class of homomorphisms $\rho:\pi_{1}(M)\longrightarrow
T_{2}^{u}$ whose diagonal elements $\rho_{ii}$ satisfy:
$\rho_{ii}=exp(\widetilde\rho_{ii})$ for some
$\widetilde\rho_{ii}\in Hom(\pi_{1}(M),\Co)$, $i=1,2$.
(E.g., if $H^{2}(M,\Z)$ is torsion free then each element of
$Hom(M,T_{2}^{u})$ belongs to $S_{2}^{u}(M)$.)
In what follows 
$f:M_{1}\longrightarrow M_{2}$ is a complex surjective mapping
of compact K\"{a}hler manifolds and $G', G''$ denote the first and the
second commutant groups of a group $G$. 
\begin{Th}\label{te5}
Assume that for any 
$\tau\in S_{2}^{u}(M_{1})$ there is 
$\tau'\in S_{2}^{u}(M_{2})$ such that
$\tau=\tau'\circ f_{*}$. Then for any $\rho\in S_{n}(M_{1})$ there exists
$\rho'\in S_{n}(M_{2})$ such that
$\rho=\rho'\circ f_{*}$.
\end{Th}
\begin{R}\label{Camp}
{\em A result similar to Theorem
\ref{te5} is also valid in the case of representations generated by 
connections from ${\cal A}_{n}^{t}$
with nilpotent (0,1)-components. 
In this case it suffices to assume that} $f$ {\em induces an
isomorphism of} $H_{1}(M_{1},\Re)$ {\em and} $H_{1}(M_{2},\Re),$ {\em
see [Br].
This assumption holds, e.g., for} $f$ {\em being a smoothing of the 
Albanese map} $\alpha_{M}$ {\em of a compact K\"{a}hler manifold} $M$ 
{\em (here} $M_{1}, M_{2}$ 
{\em be a desingularizations of} $M$ {\em and} $\alpha_{M}(M)\subset Alb(M),$
{\em 
respectively). Then the analog of Theorem \ref{te5} implies 
(see, e.g. [ABCKT], Proposition 3.33):}\\
{\bf Theorem} (Campana). $f$  induces an isomorphism of the
de Rham fundamental groups of $M_{1}$ and $M_{2}$.
\end{R}

Let us introduce now the class $S$ of compact K\"{a}hler manifolds $M$ for
which $\displaystyle \cup_{n\geq 1}S_{n}(M)$ separates elements of
$\pi_{1}(M)$. 
\begin{Th}\label{te6}
Assume that $M_{1}\!\in \!S$ and
$f$ induces an isomorphism of $\pi_{1}(\!M_{1}\!)/\!\pi_{1}(\!M_{1}\!)''$ and
$\pi_{1}(M_{2})/\pi_{1}(M_{2})''$. Then $f_{*}$ imbeds $\pi_{1}(M_{1})$ as a
subgroup of a finite index in $\pi_{1}(M_{2})$.
\end{Th}

In a forthcoming paper we will demonstrate another application of 
Theorem \ref{maintheo}.\\
{\bf Theorem.}
{\em
Assume that} $M\in S$ {\em satisfies}
$$
\begin{array}{lr}
(i)\ \ \ \ \pi_{2}(M)=0;\\
(ii)\ \ \ dim_{\Co}M\geq\frac{1}{2}rank(\pi_{1}(M)'').
\end{array}
$$
{\em Then}\\
{\em (a)}\ \ $dim_{\Co}M=\frac{1}{2}rank(\pi_{1}(M)'');$\\
{\em (b)}\ \ $\pi_{1}(M)$ {\em
is isomorphic to a lattice in a Lie group} $G$ {\em which is a semidirect
product of} $\Co^{m}$ {\em and} $\Re^{2k}$ {\em 
determined by a unitary representation} $\Re^{2k}\longrightarrow U_{m}(\Co)$.
{\em Here} $2m+2k=rank(\pi_{1}(M)'')$ {\em and}
$2m=rank(\pi_{1}(M)')$.\\
This gives, in particular, a classification of 
compact solvmanifolds admitting a K\"{a}hler structure.

At the end of the paper we will show that the 
results formulated above hold also for the class of manifolds dominated
by a compact K\"{a}hler manifold.
\sect{\hspace{-1em}. Theorem \ref{maintheo}: Principal Results.}    
{\bf 2.1.} 
Theorem \ref{maintheo} follows from
results formulated below. In order to formulate the first of them recall
that any flat connection $\omega$ on a topologically trivial complex vector 
bundle $M\times\Co^{n}$ (over a compact K\"{a}hler manifold $M$) is 
determined by equation
\begin{equation}\label{form1}
df=\omega f \ \ \ \ \ (f\in C^{\infty}(M,GL_{n}(\Co))),
\end{equation}
with $\omega$ satisfying (\ref{form2}) (the condition of local solvability). 
For a family $\{f_{i}\}_{i\in I}$ of local solutions of (\ref{form1})
defined on an open covering  $\{U_{i}\}_{i\in I}$ the flat structure on
$M\times\Co^{n}$ is determined by locally constant cocycle
$\{c_{ij}:=f_{i}^{-1}f_{j}\}_{i,j\in I}$. Further, we can 
rewrite (\ref{form1}) in the equivalent form 
\begin{equation}\label{form4}
\partial f=\omega_{1}f,
\end{equation}
\begin{equation}\label{form5}
\overline{\partial}f=\omega_{2}f,
\end{equation}
$\omega_{1}$ and $\omega_{2}$ being a (1,0)-form and a (0,1)-form, 
respectively, 
and $\omega=\omega_{1}+\omega_{2}$. As follows from (\ref{form2}) 
the system (\ref{form4})-(\ref{form5}) is locally solvable. 
It is worth pointing out that the local solvability of each of these 
equations separately is equivalent
to the fulfillment of one of the corresponding conditions:
\begin{equation}\label{form6}
\partial\omega_{1}-\omega_{1}\wedge\omega_{1}=0,
\end{equation}
\begin{equation}\label{form7}
\overline{\partial}\omega_{2}-\omega_{2}\wedge\omega_{2}=0.
\end{equation}
Our first result related to the following\\
{\em Complement problem}. Given $\omega_{2}$ satisfying (\ref{form7}), 
find $\omega_{1}$ for which the system (\ref{form4})-(\ref{form5}) is locally
solvable.
\begin{Th}\label{Te1}
Suppose that $\omega_{2}$ is a triangular (0,1)-form 
satisfying (\ref{form7}). Then there exists a triangular 
(1,0)-form $\omega_{1}$ such that $\omega=\omega_{1}+\omega_{2}\in
{\cal A}_{n}^{t}$, i.e., satisfies
(\ref{form2}).
In addition, there exists a $T_{n}(\Co)$-valued $d$-gauge transform sending 
$\omega$ to a triangular 1-form $\eta=\eta_{1}+\eta_{2}$ such that  
$$
diag(\eta_{2})=-\overline{\eta_{1}}.
$$ 
Here $diag(\phi)$ is the diagonal of $\phi$, 
and $\overline{\phi}$ denotes the complex conjugate of $\phi$.
\end{Th}
{\bf 2.2.} Let $E$ be a flat vector bundle over $M$ of complex rank $n$
constructed by a locally constant cocycle $\{c_{ij}\}_{i,j\in I}$ defined on 
an open covering $\{U_{i}\}_{i\in I}$.
Further, let $End(E)$ be the vector bundle of linear endomorphisms of $E$.
According to (\ref{section}) the operators $d$ and $\wedge$ are well-defined
on the set of matrix-valued 1-forms with values in $End(E)$. In particular,
it makes sense to consider 1-forms satisfying an equation similar to
(\ref{form2}). Let $h$ be a linear $C^{\infty}$-authomorphism of $E$
determined by a family 
$\{h_{i}\}_{i\in I}$ $(h_{i}\in C^{\infty}(U_{i}, GL_{n}(\Co)))$
satisfying
$$
h_{j}=c_{ij}^{-1}h_{i}c_{ij}\ \ \ on\ \ \ U_{i}\cap U_{j}.
$$
Then a $d$-gauge transform $d_{h}^{E}$ defined on the set
of matrix-valued 1-forms $\alpha$ with values in $End(E)$ is given
by a formula similar to (\ref{form3})
$$
d_{h}^{E}(\alpha)=h^{-1}\alpha h-h^{-1}dh.
$$
Clearly, $d_{h}^{E}$ preserves the class of 1-forms satisfying  
an $End(E)$-valued equation (\ref{form2}). Let now
$E\in {\cal U}_{\oplus}^{n}$ and $h$ belong to $Aut_{\infty}^{t}(E)$,
the group of triangular $\penalty -10000$
$C^{\infty}$-authomorphisms of $E$.
Then $d_{h}^{E}$ preserves also the class of $End(E)$-valued 1-forms
with nilpotent (0,1)-components. Since $E$ is a direct sum of 
topologically trivial vector bundles $M\times\Co$, the group
$Aut_{\infty}^{t}(E)$ is isomorphic to $C^{\infty}(M,T_{n}(\Co))$.
In what follows we identify these two groups.

We now proceed to describe the correspondence map from Theorem \ref{maintheo}.
We let ${\cal E}_{\psi}$ denote the class of connections 
from ${\cal A}_{n}^{t}$ such that diagonals
of their (0,1)-components equal $\psi$.
\begin{Proposition}\label{Prop1}
(1) For every $\overline{\partial}$-closed (0,1)-form $\psi$ there is an
invertible diagonal matrix function $h_{\psi}$ such that
$d_{h_{\psi}}({\cal E}_{\psi})={\cal E}_{\widetilde \psi}$, where
$\widetilde \psi$ is the harmonic component in the Hodge decomposition
of $\psi$.\\
(2) For every diagonal harmonic (0,1)-form $\psi$
there are a vector bundle $E_{\psi}$ over $M$ and 
an injective mapping $\tau_{\psi}$  of ${\cal E}_{\psi}$ to the set of
$End(E_{\psi})$-valued 1-forms such that\\
(a)\ \ $E_{\psi}\in {\cal U}_{\oplus}^{n}$;\\ 
(b)\ \ $\tau_{\psi}\circ d_{g}=d_{g}^{E_{\psi}}\circ\tau_{\psi}$
for every $g\in C^{\infty}(M,T_{n}(\Co))$;\\
(c)\ \ $\tau_{\psi}({\cal E}_{\psi})$ consists of forms with
nilpotent (0,1)-components satisfying (\ref{form2}).
\end{Proposition}
The proof of the proposition will also show that every element of
${\cal U}_{\oplus}^{n}$ coincides with some $E_{\psi}$.

According to the above proposition moduli space 
${\cal B}_{n}^{t}$ of flat
connections from ${\cal A}_{n}^{t}$ is isomorphic to
a similar moduli space of forms from $\tau_{\psi}({\cal E}_{\psi})$.
\begin{Proposition}\label{Tebun}
For every $\eta\in \tau_{\psi}({\cal E}_{\psi})$ there is a transform
$d_{g}^{E_{\psi}}$ with $g\in Aut_{\infty}^{t}(E_{\psi})$ such that
$$
d_{g}^{E_{\psi}}(\eta)=\widetilde\eta_{1}+\widetilde\eta_{2}, 
$$
where $\partial\widetilde\eta_{1}=0$ and $\widetilde\eta_{2}$ is 
a $d$-closed nilpotent antiholomorphic form.
\end{Proposition}

This result implies that $\widetilde\eta_{1}$ can be decomposed into the sum 
$\alpha+\partial h$, 
where $\alpha$ is its harmonic component in the Hodge decomposition. 
It is worth noting that $\widetilde\eta_{2}$ and $\alpha$ belong
to the space ${\bf H}_{d}^{1}(End(E_{\psi}))$ of $d$-harmonic forms
determined in Introduction (see Proposition \ref{Te17} below).
Moreover, condition (\ref{form2}) together with the 
$\partial\overline{\partial}$-lemma (see Lemma \ref{Te18} below)
imply that $[\alpha +\widetilde\eta_{2},\alpha +\widetilde\eta_{2}]$ 
represents 0 in the de Rham cohomology group 
$H^{2}(M,End(E_{\psi}))$. 
The converse of the latter statement is also true.
Namely, let $\alpha$ be an $End(E_{\psi})$-valued
$d$-harmonic (1,0)-form and $\theta$ be a 
$d$-harmonic nilpotent $End(E_{\psi})$-valued (0,1)-form.
\begin{Proposition}\label{Te3}
Let $[\alpha+\theta,\alpha+\theta]$ represent zero in 
$H^{2}(M, End(E_{\psi}))$.
Then there exists a unique up to a flat additive summand section $h$ such that 
$(\alpha+\partial h)+\theta$ satisfies (\ref{form2}).
\end{Proposition}

Finally, to complete Theorem \ref{maintheo} we have to prove
the following uniqueness result.
\begin{Proposition}\label{Te4}
Let $\alpha_{1},\beta_{1}$ and $\alpha_{2},\beta_{2}$ be 
$End(E_{\psi})$-valued (1,0)- and (0,1)-forms, respectively. Suppose that\\
(a)\ \ $\alpha_{1}+\alpha_{2}$ and $\beta_{1}+\beta_{2}$ belong to
$\tau_{\psi}({\cal E}_{\psi})$ and are $d$-gauge equivalent;\\
(b)\ \ $\alpha_{2},\beta_{2}$ are $d$-closed nilpotent forms;\\
Then the $d$-gauge equivalence is defined by a flat automorphism 
of $E_{\psi}$.
\end{Proposition}

In other words, $\widetilde\eta_{1}+\widetilde\eta_{2}$ in Proposition 
\ref{Tebun} is unique up to conjugation by flat automorphisms. 
We now summarize the above results.

The space ${\cal B}_{n}^{t}$ is isomorphic to disjoint union of
the moduli spaces $\penalty-10000$
${\cal E}_{\psi}/C^{\infty}(M,T_{n}(\Co))$ with
diagonal harmonic (0,1)-forms $\psi$. Further, the mapping $\tau_{\psi}$
defines an isomorphism between ${\cal E}_{\psi}/C^{\infty}(M,T_{n}(\Co))$
and $\tau_{\psi}({\cal E}_{\psi})/Aut_{\infty}^{t}(E_{\psi})$. The latter,
in turn, is isomorphic to ${\cal S}_{E_{\psi}}^{n}:=
{\bf H}_{0}^{t}(End(E_{\psi}))/Aut_{f}^{t}(E_{\psi})$. This completes
the description of the correspondence of Theorem \ref{maintheo}.
\begin{R}\label{gold}
{\rm It was proved 
by Goldman and Millson [GM] and independently by Simpson [S] that the 
representation varieties of K\"{a}hler groups have at worst quadratic 
singularities at reductive representations. Theorem \ref{maintheo} shows
that this ``quadratic law'' is also of global nature if we restrict
ourselves to some naturally determined subsets of ${\cal M}_{GL_{n}(\Co)}$.}
\end{R}
\sect{\hspace{-1em}. Auxiliary Results.}
{\bf 3.1.} Let $D$ be one of the operators $d$, $\overline{\partial}$ or 
$\partial$. If $g\subset gl_{n}(\Co)$ is the Lie algebra of a Lie group 
$G\subset GL_{n}(\Co)$ then we denote
by ${\cal A}_{D}(g)$ the space of locally 
integrable $D$-connections in the principle bundle $M\times G$ over $M$
defined by 
\begin{equation}\label{form10}
Df=\omega f \ \ \ \ \ (f\in C^{\infty}(M,G))
\end{equation}
with a $g$-valued differential forms $\omega$. 
The condition of integrabilty of a connection is
$$
D\omega-\omega\wedge\omega=0.
$$

Let ${\cal B}_{D}(g)$ denote the moduli space of ${\cal A}_{D}(g)$,
i.e., the set of $D$-gauge equivalent classes of 
connections from ${\cal A}_{D}(g)$.
Further we introduce the class ${\cal V}_{D}(G)$ of isomorphic 
$G$-topologically trivial vector bundles with $D$-trivial cocycles 
$\{c_{ij}\}$ (this means that the principle $G$-bundle constructed by this
cocycle is topologically trivial and $Dc_{ij}=0$ for all $i,j$). 
In particular, $\{c_{ij}\}$ is holomorphic for $D=\overline{\partial}$, 
locally constant for $D=d$, and antiholomorphic for $D=\partial$.

Then there exists bijection 
$$
i_{D}:{\cal B}_{D}(g)\longrightarrow {\cal V}_{D}(G)
$$ 
defined in the  following way (see, e.g., [O], sect. 5, 6 for details).
Let $\{U_{i}\}_{i\in I}$ be an open covering of $M$ and 
$f_{i}\in C^{\infty}(U_{i},G)$ be a solution of (\ref{form10}) on $U_{i}$.  
If we set $c_{ij}=f_{i}^{-1}f_{j}$ then 
$\{c_{ij}\}$ is a $D$-trivial cocycle and so it determines an element of 
${\cal V}_{D}(G)$. The construction is independent of the choice of the 
element 
of an equivalence class in ${\cal B}_{D}(g)$ and, therefore, it correctly 
defines the required mapping $i_{D}$.  
For an $\omega\in {\cal A}_{D}(g)$ we let $[\omega]\in {\cal B}_{D}(g)$ 
denote its $D$-gauge equivalence class.

Since each locally constant cocycle is holomorphic and antiholomorphc simultaneously, 
the identity mapping induces natural mappings 
\begin{equation}\label{h}
h:{\cal V}_{d}(G)\longrightarrow {\cal V}_{\overline{\partial}}(G)\ \ 
and\ \  
\overline{h}:{\cal V}_{d}(G)\longrightarrow {\cal V}_{\partial}(G).
\end{equation}
Namely, if {\bf E} is the sheaf of locally constant sections of a 
vector bundle $E\in {\cal V}_{d}(G)$ then vector bundles
$h(E)$ and $\overline{h}(E)$ are determined by
sheaves ${\bf E}\otimes_{\Co} {\cal O}_{M}$ and $\penalty -10000$
${\bf E}\otimes_{\Co}\overline{{\cal O}}_{M}$, respectively. 

It is worth noting that the moduli space of isomorphic vector 
bundles with locally constant $G$-cocycles 
({\em flat bundles}) is isomorphic to the quotient 
${\cal M}_{G}:=Hom(\pi_{1}(M),G)/G$ of the space of
representations of $\pi_{1}(M)$ in $G$, by the action of $G$ given
by conjugation (see, e.g., [KN], Ch.2, sect.9). 
\begin{Proposition}\label{Te10}
Let ${\omega_{2}}\in {\cal A}_{\overline{\partial}}(gl_{n}(\Co))$. Then the 
following statements are equivalent:\\
(i) there exists a $gl_{n}(\Co)$-valued (1,0)-form $\omega_{1}$ such that 
$\omega=\omega_{1}+\omega_{2}$ belongs to ${\cal A}_{d}(gl_{n}(\Co))$; \\
(ii) there exists an element $E\in {\cal V}_{d}(GL_{n}(\Co))$ such that 
$$
h(E)=i_{\overline{\partial}}([\omega_{2}]).
$$
\end{Proposition}
{\bf Proof.} Let $\Pi_{0,1}:{\cal E}^{1}(M)\otimes gl_{n}(\Co)\longrightarrow
{\cal E}^{0,1}(M)\otimes gl_{n}(\Co)$ be the projection from the space of
matrix-valued 1-forms defined on $M$ onto the space of (0,1)-forms 
induced by the type decomposition. Clearly, $\Pi_{0,1}$ maps  
${\cal A}_{d}(gl_{n}(\Co))$ in ${\cal A}_{\overline{\partial}}(gl_{n}(\Co))$
and commutes with actions of the corresponding gauge transform groups.
Denote by $\widetilde\Pi_{0,1}:
{\cal B}_{d}(gl_{n}(\Co))\longrightarrow
{\cal B}_{\overline{\partial}}(gl_{n}(\Co))$
the mapping induced by $\Pi_{0,1}$.
Then the required statement follows from the commutative diagram
\begin{equation}\label{diag1}
\begin{array}{ccccc}
{\cal B}_{d}(gl_{n}(\Co))&\stackrel{\widetilde\Pi_{0,1}}{\longrightarrow}
&{\cal B}_{\overline{\partial}}(gl_{n}(\Co))& &\\
i_{d}\downarrow& &\downarrow i_{\overline{\partial}}& &\\
{\cal V}_{d}(GL_{n}(\Co))&\stackrel{h}{\longrightarrow}&
{\cal V}_{\overline{\partial}}(GL_{n}(\Co))& &\Box
\end{array}
\end{equation}
{\bf 3.2.}
Below we denote by ${\cal V}$ the category of vector bundles equipped 
with one of the following structures: $C^{\infty}$, holomorphic, 
antiholomorphic or flat. If $E\in {\cal V}$ then
{\bf E} denotes the sheaf of its local sections determining the structure
of $E$.

Let now $E,\ E_{1},\ E_{2}$ belong to ${\cal V}$.
\begin{D}\label{Te11}
$E$ is said to be an extension  of $E_{2}$ by $E_{1}$ if the sequence 
\begin{equation}\label{extens}
0\longrightarrow E_{1}\longrightarrow E\longrightarrow E_{2}\longrightarrow 0
\end{equation}
is exact.\\
\newcommand{\dia}[9]
{
\begin{array}{l}
0\longrightarrow #1_{1}\longrightarrow #2\longrightarrow #3_{2}\longrightarrow 0 \\
\ \       \ \  \     #4_{1}\downarrow \ \ \ \ #5\downarrow \ \ \ \ #6_{2}\downarrow \\
0\longrightarrow #7_{1}\longrightarrow #8\longrightarrow #9_{2}\longrightarrow 0
\end{array}
}
Extensions $E$ of $E_{2}$ by $E_{1}$ and $F$ of $F_{2}$ by $F_{1}$ are
isomorphic in ${\cal V}$ if there exists a commutative diagram
\begin{equation}\label{form11}
\dia{E}{E}{E}{j}{j}{j}{F}{F}{F}
\end{equation}
where $j_{1}, j, j_{2}$ are isomorphisms of the corresponding ${\cal V}$-bundles.\\
In the case of $j_{1}=id$ and $j_{2}=id$ these extensions are called  
equivalent.
\end{D}
Let $E$ be an extension of $E_{2}$ by $E_{1}$. Then (\ref{extens}) induces
the exact sequence 
$$
0\longrightarrow Hom(E_{2},E_{1})\longrightarrow Hom(E_{2},E)
\longrightarrow Hom(E_{2},E_{2})\longrightarrow 0
$$
(here all bundles have the same structure as $E_{1}$ and $E_{2}$).
The above sequence, in turn, induces the exact sequence of 
\v{C}ech cohomology groups of the corresponding sheaves
$$
\begin{array}{l}
0\longrightarrow H^{0}(M,Hom({\bf E_{2},E_{1}}))\longrightarrow 
H^{0}(M,Hom({\bf E_{2},E}))\longrightarrow \\
H^{0}(M,Hom({\bf E_{2},E_{2}}))
\stackrel{\delta}{\longrightarrow}H^{1}(M,Hom({\bf E_{2},E_{1}}))
\longrightarrow ...
\end{array}
$$
Let $I\in H^{0}(M,Hom({\bf E_{2},E_{2}}))$ be the identity section. 
Then it is well known that $\delta (I)$ uniquely determines the 
class of extensions of $E_{2}$ by $E_{1}$ equivalent to $E$.
\begin{Proposition}\label{Te12}
([A], Proposition 2). The equivalence classes of extensions of $E_{2}$ by 
$E_{1}$ are in one-to-one correspondence with the elements of 
$H^{1}(M,Hom({\bf E_{2}},{\bf E_{1}}))$ and the trivial extension corresponds
to the trivial element.\ \ \ $\Box$ 
\end{Proposition}
\begin{R}\label{Te13}
{\em It follows directly from Definition \ref{Te11} 
that if} $E_{i}\in {\cal V}_{D}({GL_{k}}_{i}(\Co))$,
$i=1,2$, {\em then} $E\in {\cal V}_{D}(G)$, {\em 
where the structure group $G$ consists 
of elements of the form}
\[
\left(\begin{array}{cc}A_{1}& * \\ 0&A_{2}\end{array}\right)
\]
{\em with} $A_{i}\in {GL_{k}}_{i}(\Co)$, $i=1,2$.
\end{R}

Let now $E$ and $F$ be isomorphic extensions of $E_{2}$ by $E_{1}$ and 
$F_{2}$ by $F_{1}$, respectively. Let $k_{i}$ be the rank of
$E_{i}$, $i=1,2$, and $G$ be the Lie group from the above remark. Consider
principle bundles $E_{G}$ and $F_{G}$ with the structure group $G$ 
corresponding to $E$ and $F$. Then it follows immediately from the definitions
that 
\begin{Proposition}\label{Te14}
Any isomorphism $j:E\longrightarrow F$ 
determined by (\ref{form11})
induces an isomorphism $j_{G}$ of $G$-bundles $E_{G}$ and $F_{G}$. 
Moreover, restriction of $j_{G}$ to a fibre is determined as 
left multiplication by an element of $G$.\ \ \ $\Box$
\end{Proposition}

Consider now an extension $E$ of $E_{2}$ by $E_{1}$ 
in the category of flat bundles. (So structure group $G$ of $E$ is now
defined as in Remark \ref{Te13}.) In this case
the natural mappings 
$h:{\cal V}_{d}(G)\longrightarrow {\cal V}_{\overline{\partial}}(G)$ 
and $\overline{h}:{\cal V}_{d}(G)\longrightarrow {\cal V}_{\partial}(G)$, see 
({\ref{h}}), determine extensions $h(E)$ of 
$h(E_{2})$ by $h(E_{1})$
and $\overline{h}(E)$ of $\overline{h}(E_{2})$ by $\overline{h}(E_{1})$.
According to Proposition \ref{Te12} and the Dolbeault theorem the former 
extension is defined by an element of the group 
$H^{1}(M,Hom({\bf E_{2}}\otimes_{\Co}{\cal O}_{M},{\bf E_{1}}\otimes_{\Co}{\cal O}_{M}))$,
and each element of this group is given by a $\overline{\partial}$-closed
(0,1)-form with values in $Hom(E_{2},E_{1})$. The latter extension is defined  
in the same way by a $\partial$-closed (1,0)-form with values 
in $Hom(E_{2},E_{1})$.\\
The elements of the cohomology groups that appeared here should be found 
as follows.

Let $\eta\in H^{1}(M,Hom({\bf E_{2},E_{1}}))$ be an element defining the 
extension $E$. Let $\Pi_{0,1}$, $\Pi_{1,0}$ be the natural projections 
from the space of 1-forms onto spaces of (0,1)- and (1,0)-forms, respectively.
By the same symbols we denote mappings of the corresponding cohomology 
groups induced by $\Pi_{0,1}$ and $\Pi_{1,0}$. So that
$$
\begin{array}{l}
\Pi_{0,1}(\eta)\in H^{1}(M, Hom({\bf E_{2}}\otimes_{\Co}{\cal O}_{M},
{\bf E_{1}}\otimes_{\Co}{\cal O}_{M})),\\
\Pi_{1,0}(\eta)\in 
H^{1}(M, Hom({\bf E_{2}}\otimes_{\Co}\overline{{\cal O}}_{M},
{\bf E_{1}}\otimes_{\Co}\overline{{\cal O}}_{M})).
\end{array}
$$
\begin{Proposition}\label{Te15}
The classes of extensions equivalent to $h(E)$ and $\overline{h}(E)$
are uniquely defined by $\Pi_{0,1}(\eta)$ and $\Pi_{1,0}(\eta)$, respectively.
\end{Proposition}
{\bf Proof.} In the case of $h(E)$ the result follows directly from
de Rham's and Dolbeault's theorems applied to the second 
column of the commutative diagram
\[
\begin{array}{ccc}
H^{0}(M,Hom({\bf E_{2},E_{2}}))&
\stackrel{\delta}{\longrightarrow}&H^{1}(M,Hom({\bf E_{2},E_{1}}))\\
h\downarrow & &\downarrow\Pi_{0,1}\\
H^{0}(M,Hom({\bf E_{2}}\otimes_{\Co}{\cal O}_{M},
{\bf E_{2}}\otimes_{\Co}{\cal O}_{M}))&
\stackrel{\delta}{\longrightarrow}&
H^{1}(M,Hom({\bf E_{2}\otimes_{\Co}{\cal O}_{M}},
{\bf E_{1}}\otimes_{\Co}{\cal O}_{M})).
\end{array}
\]
The case of $\overline{h}(E)$ is similar.\ \ \ \ \ $\Box$
\\
{\bf 3.3.} In this part we collect several facts on the class
${\cal S}{\cal B}$ of bundles with  
connected solvable complex Lie groups as structure groups.\\
(a) ${\cal S}{\cal B}$ is closed 
under tensor products and duality, i.e., $E^{*}$ and $E\otimes D$ belong to 
${\cal S}{\cal B}$ together with $E,D$.\\
(b) Every element $E\in {\cal S}{\cal B}$ can be thought of as a vector bundle 
with structure group $T_{n}(\Co)$ (for some $n$).\\
Actually, according to the Lie theorem, 
for any connected solvable 
subgroup $G$ of $GL_{n}(\Co)$ there exists a matrix $B\in GL_{n}(\Co)$ such 
that 
$B^{-1}GB$ is imbedded as a subgroup in the group $T_{n}(\Co)$. Moreover, 
let $E$ have one of the structures: holomorphic, 
antiholomorphic, or flat. Then the above
transform generates an isomorphism of $E$ preserving this structure.\\
(c) Every $E\in {\cal S}{\cal B}$ is the result of successive extensions of 
bundles with triangular structure groups by means of rank 1 vector 
bundles.\\ 
Indeed, for the action of $T_{n}(\Co)$ on $\Co^{n}$ there exists
a one-dimensional
invariant subspace such that $T_{n-1}(\Co)$ acts on the factor 
space. Therefore $E$ with structure group  $T_{n}(\Co)$
is an extension of the bundle $E_{n-1}$ by the bundle $E_{1}$; here $E_{i}$ 
has structure group $T_{i}(\Co)$, $i=1,\ n-1$. 

Let 
$$
\begin{array}{c}
\{0\}=E_{0}\subset E_{1}\subset E_{2}\subset ... \subset E_{n-1}\subset E_{n}=E,
\\
\{0\}=F_{0}\subset F_{1}\subset F_{2}\subset ... \subset F_{n-1}\subset F_{n}=F
\end{array}
$$ 
be isomorphic flags of bundles with triangular structure groups. According to
Proposition \ref{Te14} the above isomorphism is defined (in corresponding local
coordinates on $E$ and $F$) by triangular matrices.

Let now
$$
Gr^{*} E:=\bigoplus_{i=1}^{n} E_{i}/E_{i-1}
$$
be the associated graded vector bundle with cocycle 
defined as the diagonal of the cocycle of $E$.\\
(d) $E$ is isomorphic to $Gr^{*} E$ in the category of $C^{\infty}$-bundles.\\
Really, by Proposition \ref{Te12} every vector bundle $E$ over $M$ with 
structure group $T_{n}(\Co)$ is defined by $E_{1},\ E_{n-1}$
and an element $H^{1}(M,Hom({\bf E_{n-1}},{\bf E_{1}}))$. But the latter group
is trivial in the category of $C^{\infty}$-bundles, because 
$Hom({\bf E_{n-1}},{\bf E_{1}})$ is a fine sheaf. 

As a corollary we have the following statement\\
(e) Every bundle $E\in {\cal V}_{D}(T_{n}(\Co))$ is $T_{n}(\Co)$-isomorphic to 
the direct sum of topologically trivial vector bundles $M\times\Co$.\\
(f) The class $\displaystyle \cup_{n\geq 1}{\cal V}_{D}(T_{n}(\Co))$ is  
closed under tensor products and duality. 
\\
{\bf 3.4.} In this part we recall some facts of Hodge 
theory.

Let $E$ be a flat vector bundle with structure group $U_{n}(\Co)$
over a compact K\"{a}hler manifold $M$. Then the operator of differentiation
$d$ is well-defined on the set ${\cal E}(E)$ of $E$-valued forms and
determines a connection on $E$ compatible with the complex structure and the 
flat Hermitian metric on $E$. 
Let $Z_{d}^{p,q}(E)$ be the space of $d$-closed $E$-valued $(p,q)$-forms.
As usual, one defines the cohomolgy groups of $E$ by
$$
\begin{array}{c}
H^{p,q}(E):=Z_{d}^{p,q}(E)/(d{\cal E}(E)\cap Z_{d}^{p,q}(E)), \ \ \ 
{\bf H}^{p,q}(E):=\{\eta\in {\cal E}^{p,q}(E)\;\ \ \Delta_{d}\eta=0\}, \\ 
{\bf H}_{d}^{r}:=\{\eta\in {\cal E}^{r}(E)\;\ \ \Delta_{d}\eta=0\},
\end{array}
$$
where $\Delta_{d}$ denotes the $d$-Laplacian on $E$.\\ 
Let $H^{r}(M,{\bf E})$ denote the \v{C}ech cohomology of the sheaf 
{\bf E} of locally constant sections of $E$.
\begin{Proposition}\label{Te17}
{\em (the Hodge decomposition)} \\ 
$$
\begin{array}{c}
\displaystyle H^{r}(M,{\bf E})\cong \bigoplus_{p+q=r}{\bf H}^{p,q}(E)\cong 
\bigoplus_{p+q=r}H^{p,q}(E), \ \ \ \ 
\overline{{\bf H}^{p,q}(E)}\cong {\bf H}^{q,p}(E^{*}). 
\end{array}
$$
\end{Proposition}
The proof follows from K\"{a}hler's identities for the connection $d$,
see, e.g., [ABCKT], p. 104, which give the identities between Laplacians
$$
\Delta_{d}=2\Delta_{\partial}=2\Delta_{\overline{\partial}}
$$ 
where $\Delta_{\partial}$ and $\Delta_{\overline{\partial}}$ are $\partial$- 
and $\overline{\partial}$- Laplacians on $E$.

These identities and the Dolbeault theorem give also the isomorphisms
$$
{\bf H}^{p,q}(E)\cong H_{\overline{\partial}}^{p,q}(E)\cong H^{q}(M,\Omega^{p}_{M}\otimes_{\Co}{\bf E})
$$
where $\Omega^{p}_{M}$ is the sheaf of germs of holomorphic 
$p$-forms on $M$.

Arguing as in the proof of the lemma in sect. 2 of Ch.1 of [GH] and applying
the very same identities we obtain
\begin{Lm}\label{Te18}
{\em ($\partial\overline{\partial}$-lemma)} Let $E$ be a flat bundle with 
structure group $U_{n}(\Co)$. Suppose that $\omega$ is a $d$-closed 
$E$-valued $(p,q)$-form which is $\partial$- or $\overline{\partial}$- exact. 
Then there exists an $E$-valued $(p-1,q-1)$-form $\kappa$ such that 
$$
\omega=\partial\overline{\partial}(\kappa).
$$
\end{Lm}
{\bf 3.5.} In this part we collect several facts on relations between
equations of type (\ref{form1}) and vector bundles
$Hom(E_{1},E_{2})$. 

We begin with equation
\begin{equation}\label{form12}
df=\omega_{1}f-f\omega_{2}
\end{equation}
with $\omega_{1},\omega_{2}$ satisfying (\ref{form2}).
The right side can be written as
$(1\otimes \omega_{1}-\omega_{2}^{t}\otimes 1)f$, where $f$ is now thought of
as $n^{2}$-vector. The mapping $i_{d}$ in the following proposition is
defined as in Section 3.1.
\begin{Proposition}\label{Te19}
$i_{d}(1\otimes \omega_{1}-\omega_{2}^{t}\otimes 1)$ is a flat vector bundle 
isomorphic to\\ $Hom(i_{d}(\omega_{2}),i_{d}(\omega_{1}))$.
\end{Proposition}
{\bf Proof}. Let $\{U_{i}\}_{i\in I}$ be an open covering of $M$ and 
$f_{ki}\in C^{\infty}(U_{i},GL_{n}(\Co))$ be a solution on $U_{i}$ of equation 
(\ref{form1}) with $\omega=\omega_{k}$ $(k=1,2)$. Then 
$$
\begin{array}{l}
d((f_{2i}^{t})^{-1}\!\otimes\! f_{1i})=-\omega_{2}^{t}(f_{2i}^{t})^{-1}\!\otimes\! f_{1i}+
(f_{2i}^{t})^{-1}\!\otimes\!\omega_{1}f_{1i}=(1\!\otimes\!\omega_{1}-\omega_{2}^{t}\!\otimes\! 1)
((f_{2i}^{t})^{-1}\!\otimes\! f_{1i}).
\end{array}
$$
This means that equation (\ref{form12}) is locally solvable and defines a
flat vector bundle with cocycle 
$$
\{((f_{2i}^{t})^{-1}\otimes f_{1i})^{-1}\cdot ((f_{2j}^{t})^{-1}\otimes f_{1j})\}
:=\{(c_{2ij}^{t})^{-1}\otimes c_{1ij}\}.
$$ 
Here $\{c_{kij}:=f_{ki}^{-1}f_{kj}\}$ is a cocycle determining
flat vector bundle 
$i_{d}(\omega_{k})$ $k=1,2$. Moreover,
$\{(c_{2ij}^{t})^{-1}\}$ is a cocycle determining conjugate vector 
bundle $(i_{d}(\omega_{2}))^{*}$ (see [GH], Ch.0).
This implies that $i_{d}(1\otimes\omega_{1}-\omega_{2}^{t}\otimes 1)$ is a 
flat vector bundle isomorphic to 
$(i_{d}(\omega_{2}))^{*}\otimes (i_{d}(\omega_{1}))$. But the latter 
is isomorphic to 
$Hom(i_{d}(\omega_{2}),i_{d}(\omega_{1}))$. \ \ \ \ \ $\Box$

Let now $\eta$ be a vector-valued $(p,q)$-form on $M$ satisfying 
\begin{equation}\label{form13}
\overline{\partial}\eta=\Pi_{0,1}(\omega)\wedge\eta,
\end{equation}
where $\omega$ satisfies (\ref{form2}) and
$\Pi_{0,1}$ is the natural projection from ${\cal E}^{1}(M)$ onto 
${\cal E}^{0,1}(M)$. 
Let us check that $\eta$ is a $\overline{\partial}$-closed 
$i_{d}(\omega)$-valued $(p,q)$-form.
Clearly, $\eta$ is a section of $i_{d}(\omega)$ which is 
$C^{\infty}$-isomorphic to the vector bundle $M\times\Co^{n}$ 
(for some $n$).
Further, in flat coordinates on $i_{d}(\omega)$ determined by flat
connection $\omega$, the section $\eta$ is given by the family
$$
\{\eta_{i}:=f_{i}^{-1}\eta\}_{i\in I}.
$$
Here $f_{i}$ is a local solution on $U_{i}$ of equation (\ref{form1}) 
with the form $\omega$.

From the definition of $f_{i}$ it follows that
$$
\overline{\partial}(f_{i}^{-1}\eta)=-(f_{i}^{-1}\Pi_{0,1}(\omega))\wedge\eta+
f_{i}^{-1}(\Pi_{0,1}(\omega)\wedge\eta)=0.
$$
So $\eta$ is $\overline{\partial}$-closed.\\
Applying the very same arguments in reverse order, one deduces
that each $\overline{\partial}$-closed $i_{d}(\omega)$-valued $(p,q)$-form 
given by a family
$\{\eta_{i}\}_{i\in I}$ defines a global form $\eta$ on $M$, equal to 
$f_{i}\eta_{i}$ on $U_{i}$, satisfying (\ref{form13}).\\
In the same way we can also examine the equation
\begin{equation}\label{form14}
\partial\eta=\Pi_{1,0}(\omega)\wedge\eta
\end{equation}
and prove that $\eta$ is a $\partial$-closed $i_{d}(\omega)$-valued 
$(p,q)$-form.
Here 
$\Pi_{1,0}: {\cal E}^{1}(M)\longrightarrow {\cal E}^{1,0}(M)$ is the natural 
projection.

Finally, let us consider equations 
\begin{equation}\label{form15}
\overline{\partial}\eta=\Pi_{0,1}(\omega_{1})\wedge\eta+
(-1)^{p+q+1}\eta\wedge\Pi_{0,1}(\omega_{2})
\end{equation}
\begin{equation}\label{form16}
\partial\psi=\Pi_{1,0}(\omega_{1})\wedge\psi+(-1)^{p+q+1}\psi
\wedge\Pi_{1,0}(\omega_{2})
\end{equation}
with matrix $(p,q)$-forms $\eta$ and $\psi$. They can be written in 
equivalent forms as
$$
\begin{array}{c}
\overline{\partial}\eta=(1\otimes\Pi_{0,1}(\omega_{1})-
\Pi_{0,1}(\omega_{2}^{t})\otimes 1)\wedge\eta 
\ \ \ \ \ and\\ \ \
\partial\psi=(1\otimes\Pi_{1,0}(\omega_{1})-
\Pi_{1,0}(\omega_{2}^{t})\otimes 1)\wedge\psi,
\end{array}
$$
where $\eta$ and $\psi$ are thought of 
as {\em vector} $(p,q)$-forms.\\ 
Bringing together the results proved above for such equations, one 
gets
\begin{Proposition}\label{Te20}
There exists a one-to-one correspondence between solutions of equations  
(\ref{form15}) (or (\ref{form16})) with $\omega_{i}$ satisfying 
condition (\ref{form2}) $\mbox{(i=1,\ 2)}$ and 
$\overline{\partial}$-closed ($\partial$-closed, respectively) 
$(p,q)$-forms with values in 
$Hom(i_{d}(\omega_{2}),i_{d}(\omega_{1}))$.
\ \ \ \ \ $\Box$
\end{Proposition}
\sect{\hspace{-1em}. Proof of Theorem \ref{Te1}.}
The proof is based on  Lemmas \ref{Te21} and \ref{Te22}. To formulate the
first of the results we let
$T_{n}^{u}$ denote the subgroup of elements $A\in T_{n}(\Co)$ such that
all of its diagonal elements belong to 
$U_{1}(\Co)\!\!:=\!\!\{z;\ |z|=1\}$. 
One considers a class ${\cal U}_{n}$ of flat vector bundles $F$ with
the structure group $T_{n}^{u}$ satisfying 
\begin{equation}\label{form17}
\begin{array}{l}
 
\overline{h}(F) \ \ is \ \ isomorphic\ \ to\ \ Gr^{*}\overline{h}(F)\ \ 
in \ \ the\ \  category\ \ of\\
antiholomorphic\ vector\ bundles\ with\ structure\ group\ T_{n}(\Co).
\end{array}
\end{equation}
Let ${\cal U}=\cup_{n\geq 1}{\cal U}_{n}$. Clearly, ${\cal U}$ is closed 
under tensor products and duality.

As we explained in Section 3.3 (c) any bundle $F\in {\cal U}_{n}$ is a result 
of successive extensions of flat bundles $F_{i}$ 
with structure group $T_{i}^{u}$ by flat bundles $F^{i}$ of complex rank 1 
with structure group $U_{1}(\Co)$ $(i=1,...,n)$, so that $F=F_{n}$. 
From property (\ref{form17}) it follows
that $\overline{h}(F_{i})$ is the trivial extension of $\overline{h}(F_{i-1})$
by $\overline{h}(F^{i})$. Hence the short exact sequence of sheaves of germs of
antiholomorphic $p$-forms $(p\geq 0)$ with values in the corresponding bundles
\begin{equation}\label{form18}
0\longrightarrow\overline{\Omega}^{p}(\overline{h}(F^{i}))\stackrel{\lambda}{\longrightarrow}
\overline{\Omega}^{p}(\overline{h}(F_{i}))\stackrel{\kappa}{\longrightarrow}
\overline{\Omega}^{p}(\overline{h}(F_{i-1}))\longrightarrow 0
\end{equation}
is split. 

For a flat vector bundle $F$ we let $\overline{\Omega}_{d}^{1}(F)$ denote the 
space of $F$-valued $d$-closed antiholomorphic 1-forms. The space defines a
subgroup 
[$\overline{\Omega}_{d}^{1}(F)$] of $H^{1}(M,{\bf F})$. Here {\bf F} is the
sheaf of locally constant sections of $F$. 
Further, let $\Pi_{0,1}:H^{1}(M,{\bf F})\longrightarrow
H^{1}(M,{\cal O}_{M}\otimes_{\Co}{\bf F})$ be the mapping induced by the 
projection sending a 1-form to its (0,1)-component.
\begin{Lm}\label{Te21}
Let $F\in {\cal U}$. Then the following statements hold:\\
(a) $\Pi_{0,1}:{\em [}\overline{\Omega}_{d}^{1}(F) {\em ]}\longrightarrow
H^{1}(M,{\cal O}_{M}\otimes_{\Co}{\bf F})$ is a surjection;\\
(b) every holomorphic $F$-valued $q$-form $\alpha$ is $d$-closed. If, in 
addition, $\alpha$ is $\partial$-exact, then $\alpha=0$.
\end{Lm}
{\bf Proof}. We will prove the lemma by induction on the dimension $i$ of a 
fibre of $F$.\\
(a). In case $i=1$ the structure group of $F_{1}$ is $U_{1}(\Co)$. 
Then according to the Hodge decomposition (see Section 3.4), there exists 
an isomorphism $\penalty-10000$
$f:H^{1}(M,{\cal O}_{M}\otimes_{\Co}{\bf F_{1}})\longrightarrow 
[\overline{\Omega}_{d}^{1}(F_{1})]$ 
such that $\Pi_{0,1}\circ f=id$.

Assume now that statement (a) holds for $i-1\geq 1$; we will prove it for $i$.
The definition of extensions of bundles leads to the following commutative 
diagram:
\[
\begin{array}{ccccccc}
\!\!\!\!H^{1}\!(\!M,{\cal O}_{M}\!\otimes_{\Co}\!{\bf F^{i}}\!)&
\!\!\!\!\!\!\!\!\!\!\stackrel{\lambda}{\longrightarrow}\!\!\!\!\!\!\!&
\!\!\!\!\!H^{1}\!(\!M,{\cal O}_{M}\!\otimes_{\Co}\!{\bf F_{i}}\!)&
\!\!\!\!\!\!\!\!\!\!\stackrel{\kappa}{\longrightarrow}\!\!\!\!\!\!\!&
\!\!\!\!\!H^{1}\!(\!M,{\cal O}_{M}\!\otimes_{\Co}\!{\bf F_{i-1}}\!)&
\!\!\!\!\!\!\!\!\!\!\stackrel{\delta}{\longrightarrow}\!\!\!\!&
\!\!\!\!\!\!H^{2}\!(\!M,{\cal O}_{M}\!\otimes_{\Co}\!{\bf F^{i}}\!)
\\
\Pi_{0,1}\uparrow& &
\Pi_{0,1}\uparrow& &
\Pi_{0,1}\uparrow  \\        
H^{1}(M,{\bf F^{i}})&
\stackrel{\lambda}{\longrightarrow}&H^{1}(M,{\bf F_{i}})&
\stackrel{\kappa}{\longrightarrow}&H^{1}(M,{\bf F_{i-1}})&
\stackrel{\delta}{\longrightarrow}&H^{2}(M,{\bf F^{i}})
\end{array}
\]
By de Rham's and Dolbeault's theorems each of the elements of 
these cohomology groups is represented by an $F$-valued form.
Let $\alpha$ be a $F_{i}$-valued $\overline{\partial}$-closed 
(0,1)-form representing an element of
$H^{1}(M,{\cal O}_{M}\otimes_{\Co}{\bf F_{i}})
\cong H_{\overline{\partial}}^{0,1}(M,F_{i})$. 
According to the above
diagram and the inductive hypothesis there exists a $C^{\infty}$-section $g$ 
of $F_{i-1}$ such that 
$$
\kappa(\alpha)+\overline{\partial}(g)\in \overline{\Omega}_{d}^{1}(F_{i-1}).
$$
Since $F_{i}$ is a trivial extension in the category of $C^{\infty}$-bundles,
we can find a $C^{\infty}$-section $t$ of $F_{i}$ such 
that $\kappa(t)=g$. Then $\omega:=\kappa(\alpha-\overline{\partial}t)$ is 
a $d$-closed 1-form and therefore the (1,1)-form 
$$
\alpha':=d(\alpha-\overline{\partial}t)=\partial(\alpha-\overline{\partial}t)
$$ 
can be considered 
as an $F^{i}$-valued one. Since $\lambda(\alpha')$ represents 0 in 
$\penalty -10000$ $H_{\partial}^{1,1}(M,F_{i})\cong 
H^{1}(M,\overline{\Omega}^{1}(\overline{h}(F_{i})))$ and the mapping 
$$
\lambda: H^{1}(M,\overline{\Omega}^{1}(\overline{h}(F^{i})))\longrightarrow 
H^{1}(M,\overline{\Omega}^{1}(\overline{h}(F_{i})))
$$ 
is an injection (by (\ref{form18})),  we can deduce that
$$
[\alpha']=0\in H^{1}(M,\overline{\Omega}^{1}(F^{i})).
$$ 
So $\alpha'$ 
is a $d$-closed $\partial$-exact (1,1)-form with values in a flat vector 
bundle with structure group $U_{1}(\Co)$. Then
according to the $\partial\overline{\partial}$-lemma of Section 3.4 there 
exists a $C^{\infty}$-section $s$ of $F^{i}$ such that 
$$
\partial\overline{\partial}(s)=\alpha'.
$$
We set now
$$
\beta:=\alpha-\overline{\partial}t-\overline{\partial}(\lambda(s)).
$$
Then $\beta$ is a $\overline{\partial}$-closed (0,1)-form such that 
$$
\begin{array}{c}
[\beta]=[\alpha]\in H^{1}(M,{\cal O}_{M}\otimes_{\Co}{\bf F_{i}}) \ \ \ \ \
and \ \ d\beta=0.
\end{array}
$$
Hence $\beta$ represents an element $\widetilde{\beta}$ of 
$[\overline{\Omega}_{d}^{1}(F_{i})]$
such that $\Pi_{0,1}(\widetilde{\beta})=[\alpha]$.
The proof of part $(a)$ is complete.\\
(b). We again make use of induction on $i$.  
So let $\omega_{i}$ be a $F_{i}$-valued holomorphic $q$-form. 
In case $i=1$ the Hodge identity for Laplacians (see Section 3.4) acquires 
the form
$$
\bigtriangleup_{d}(\omega_{1})=2\bigtriangleup_{\overline{\partial}}(\omega_{1})=0,
$$
and from this it follows that $d\omega_{1}=0$. If, in addition, $\omega_{1}$ is
$\partial$-exact then its $\penalty -10000$ 
$d$-harmonicity implies $\omega_{1}=0$.\\
Assume now that statement (b) holds for $i-1\geq 1$; we will prove it for $i$.
By the induction hypothesis we have $d(\kappa(\omega_{i}))=0$. But 
$d(\kappa(\omega_{i}))=\kappa(d\omega_{i})$ and therefore $d\omega_{i}$ can
be regarded as a $d$-closed 
$F^{i}$-valued holomorphic $(q+1)$-form. 

It is clear, as well, that  
$$
[d\omega_{i}]=[\partial\omega_{i}]\in H_{\partial}^{q+1,0}(M,F^{i})\cong
H^{q+1}(M,\overline{{\cal O}}_{M}\otimes_{\Co}{\bf F^{i}}).
$$
Now on account of (\ref{form18}) the mapping 
$$
\lambda:H^{q+1}(M,\overline{{\cal O}}_{M}\otimes_{\Co}{\bf F^{i}})
\longrightarrow 
H^{q+1}(M,\overline{{\cal O}}_{M}\otimes_{\Co}{\bf F_{i}})
$$ 
is an injection. On the other hand, 
$\lambda([\partial\omega_{i}])=0$
and, hence, $[\partial\omega_{i}]=0$ in $\penalty-10000$
$H^{q+1}(M,\overline{{\cal O}}_{M}\otimes_{\Co}{\bf F^{i}})$.
Taking further into account the above mentioned identity for Laplacians in 
the one-dimensional case we can deduce that $\partial\omega_{i}$ 
is a $d$-harmonic $F^{i}$-valued form. Since it is $\partial$-exact, we have
$\partial\omega_{i}=d\omega_{i}=0$.

It remains to prove that if, in addition, $\omega_{i}$ is a $\partial$-exact 
form, then it equals 0. But in this case $\kappa(\omega_{i})$ 
is a $\partial$-exact $F_{i-1}$-valued holomorphic form and, consequently,
$\kappa(\omega_{i})=0$ by the induction hypothesis. So, $\omega_{i}$ can be 
regarded as a $F^{i}$-valued holomorphic form. Moreover, according to the 
equality, 
$$
H^{q}(M,\overline{{\cal O}}_{M}\otimes_{\Co}{\bf F_{i}})=
H^{q}(M,\overline{{\cal O}}_{M}\otimes_{\Co}{\bf F^{i}})\bigoplus 
H^{q}(M,\overline{{\cal O}}_{M}\otimes_{\Co}{\bf F_{i-1}}),
$$
(see (\ref{form18})) $\omega_{i}$ is $\partial$-exact.
Thus $\omega_{i}$ is a $F^{i}$-valued $\partial$-exact holomorphic form and 
therefore
it equals 0 as we have already shown at the first step of the induction. \ \ \ \ \ 
$\Box$

Let us suppose now that $\omega_{2}$ is a triangular (0,1)-form of the class
${\cal A}_{\overline{\partial}}(t_{n})$ (see Section 3.1 for the definition 
of this class). Here $t_{n}$ denotes the Lie algebra of $T_{n}(\Co)$.
\begin{Lm}\label{Te22}
The following conditions are equivalent: \\
(i) \ For $\omega_{2}$ Theorem \ref{Te1} holds; \\
(ii)  there exists a $T_{n}^{u}$-topologically trivial flat vector bundle 
$F\in {\cal U}$ such that
$$
h(F)=i_{\overline{\partial}}(\omega_{2})(\in
{\cal V}_{\overline{\partial}}(T_{n}(\Co)) ).
$$
\end{Lm}
{\bf Proof}. $(i)\Rightarrow (ii)$. According to Theorem \ref{Te1} 
there exists a form $\eta\in {\cal A}_{d}(t_{n})$ with the canonical 
decomposition $\eta=\eta_{1}+\eta_{2}$ such that 
$$
[\omega_{2}]=[\eta_{2}]\in {\cal B}_{\overline{\partial}}(t_{n})\ \ and\ \
diag(\eta_{2})=-\overline{\eta_{1}}.
$$
Since $diag(\eta)=\eta_{1}-\overline{\eta_{1}}$ is 
$(\sqrt{-1}\cdot \Re)^{n}$-valued, the
form $\eta$ defines a unique element of ${\cal B}_{d}(t_{n}^{u})$. Here 
$t_{n}^{u}$ is the 
Lie algebra of $T_{n}^{u}$ which 
clearly consists of elements $A\in T_{n}(\Co)$
with $diag(A)\in (\sqrt{-1}\cdot \Re)^{n}$.
Therefore the flat bundle $i_{d}(\eta)$ has structure group
$T_{n}^{u}$ (see Section 3.1).

Now we make use of the identities 
$$
\overline{h}(i_{d}(\eta))=i_{\partial}(\eta_{1})\in 
{\cal V}_{\partial}(T_{n}(\Co)),\ \ \ \ h(i_{d}(\eta))=
i_{\overline{\partial}}(\omega_{2})\in 
{\cal V}_{\overline{\partial}}(T_{n}(\Co)),
$$
see Proposition \ref{Te10} for details.
But $\eta_{1}$ is a diagonal matrix form and thus the first identity 
implies that $\overline{h}(i_{d}(\eta))$ is isomorphic to 
$Gr^{*}\overline{h}(i_{d}(\eta))$ in the category of antiholomorphic vector 
bundles with structure group $T_{n}(\Co)$. 
Therefore $i_{d}(\eta)$ belongs to the class ${\cal U}$ of flat vector bundles 
with structure groups $T_{n}^{u}$ and is 
$T_{n}(\Co)$-topologically trivial 
by the definition of the class ${\cal V}_{D}(T_{n}(\Co))$. Moreover, every  
$T_{n}(\Co)$-topologically trivial vector bundle with structure group
$T_{n}^{u}$ is $T_{n}^{u}$-topologically trivial. Bearing in mind 
the second identity we deduce now that $i_{d}(\eta)$ can be taken as the bundle 
$F$ of statement $(ii)$. \\
$(ii)\Rightarrow (i)$. Let $F$ be the vector bundle of statement
$(ii)$. According to the results of Section 3.1, there exists a form 
$\theta\in {\cal A}_{d}(t_{n}^{u})$ with the canonical decomposition 
$\theta=\theta_{1}+\theta_{2}$ such that 
$$
i_{d}(\theta)=F.
$$ 
In particular, we have $diag(\theta_{1})=-diag(\overline{\theta_{2}})$. 
Moreover, as has been established in the first part of the proof 
the following equalities 
$$
i_{\partial}(\theta_{1})=\overline{h}(i_{d}(\theta))=
\overline{h}(F)=\overline{h}(Gr^{*} F)=
\overline{h}(i_{d}(diag(\theta)))=i_{\partial}(diag(\theta_{1})).
$$
hold in the class ${\cal V}_{\partial}(T_{n}(\Co))$.
This implies the existence of a $\partial$-gauge transform $\partial_{g}$ 
with a triangular matrix function $g$ such that 
$$
\partial_{g}(\theta_{1})=diag(\theta_{1}).
$$ 
Then we have for  $\psi:=d_{g}(\theta)$  the equality
$$
i_{d}(\psi)=i_{d}(\theta)=F
$$
and the first component $\psi_{1}$ of the canonical decomposition 
$\psi=\psi_{1}+\psi_{2}$
equals $\partial_{g}(\theta_{1})$, i.e., is a diagonal (1,0)-form.  
Moreover, 
$$
i_{d}(diag(\psi))=i_{d}(diag(\theta))
$$ 
in the category of flat vector bundles with the diagonal matrix structure 
group. This implies the existence of a $d$-gauge transform $d_{h}$ 
with a diagonal matrix function $h$ such that 
$$
\overline{\partial}_{h}(diag(\psi_{2}))=-\overline{\psi_{1}}.
$$
Putting now
$$
\eta:=d_{h}(\psi)
$$ 
we have defined a $t_{n}$-valued 1-form such that 
$diag(\eta_{2})=-\overline{\eta_{1}}$.  So $\eta$ satisfies
the conditions of Theorem \ref{Te1}.\\ 
It remains to define a triangular form 
$\omega$ with the second component $\omega_{2}$ in its canonical 
decomposition satisfying $\eta=d_{q}(\omega)$ for some $T_{n}(\Co)$-valued 
function $q$.
To accomplish this we note that
$$
i_{\overline{\partial}}(\theta_{2})=h(F)=i_{\overline{\partial}}(\omega_{2}),
$$
and therefore $\overline{\partial}_{p}(\theta_{2})=\omega_{2}$ for some 
$T_{n}(\Co)$-valued function $p$. If we set $\omega:=d_{p}(\theta)$
then $\omega$ satisfies the condition (\ref{form2}) because
$\theta\in {\cal A}_{d}(t_{n}^{u})$. Moreover $d_{q}(\omega)=\eta$ where 
$q:=hgp^{-1}$.\ \ \ \ \ $\Box$
\\  \\
{\bf Proof of Theorem \ref{Te1}}. Let 
$\omega_{2}\in {\cal A}_{\overline{\partial}}(t_{n})$.
According to Lemma \ref{Te22}, we have to find a $T_{n}^{u}$-topologically 
trivial flat vector bundle $F\in {\cal U}$ such that 
$$
h(F)=i_{\overline{\partial}}(\omega_{2})\in 
{\cal V}_{\overline{\partial}}(T_{n}(\Co)). 
$$
We will prove this by induction on the rank $n$ of the
holomorphic vector bundle $i_{\overline{\partial}}(\omega_{2})$.
This bundle is a result of successive extensions of holomorphic vector 
bundles
$V_{i}\in {\cal V}_{\overline{\partial}}(T_{i}(\Co))$ by rank 1
holomorphic vector bundles 
$V^{i}\in {\cal V}_{\overline{\partial}}(\Co^{*})$ $(i=1,...,n-1)$ (see 
Section 3.3). In particular, $i_{\overline{\partial}}(\omega_{2})$ is an 
extension of $V_{n-1}$ by $V^{n-1}$. 

We begin with the observation that every rank 1 holomorphic 
vector bundle 
$V\in {\cal V}_{\overline{\partial}}(\Co^{*})$ is determined by an equation 
$\overline{\partial}f=\kappa f$ with an 1-form $\kappa$ satisfying the 
condition $\overline{\partial}\kappa=0$. 
Moreover, 
a $\overline{\partial}$-gauge transform $\overline{\partial}_{g}$ in this 
case has the form
$$
\omega\mapsto \omega-g^{-1}\overline{\partial}g, \ \ \
g\in \Co^{\infty}(M,\Co^{*}).
$$
Now we are in a position to prove the result for the 1-dimensional case. Let
$V$ and $\kappa$ be as above.
Since $M$ is\ a compact K\"{a}hler  manifold, there exists a function 
$r\in C^{\infty}(M)$ 
such that $\gamma=\kappa-\overline{\partial}r$ 
is a harmonic form and, in particular, is $d$-closed. 
It is clear that $\overline{\partial}_{g}(\gamma)=\kappa$, where $g=exp(-r)$.
Let us consider now the locally solvable equation 
$$
df=(\gamma-\overline{\gamma})f.
$$
Then $d_{g}(\gamma-\overline{\gamma})=\sigma+\kappa$,
where $\sigma=-\overline{\gamma}+g^{-1}\partial g$. 
Hence, we obtain 
$$
h(i_{d}(\gamma-\overline{\gamma}))=i_{\overline{\partial}}(\kappa)=V.
$$
But $\gamma-\overline{\gamma}\in \sqrt{-1}\cdot \Re$ and therefore 
$i_{d}(\gamma-\overline{\gamma})\in {\cal V}_{d}(U_{1}(\Co)).$
It remains to set 
$$
F:=i_{d}(\gamma-\overline{\gamma}).
$$

Let us assume that the result holds for rank $n-1$; we will prove it for $n$.\\
So let $i_{\overline{\partial}}(\omega_{2})$ be an extension of $V_{n-1}$
by $V^{n-1}$. According to the induction hypothesis there exist bundles 
$F_{n-1}\in {\cal V}_{d}(t_{n-1}^{u})\cap {\cal U}$ and 
$F^{n-1}\in {\cal V}_{d}(U_{1}(\Co))$
such that 
$$
h(F_{n-1})=V_{n-1}\ \ \ and\ \ \ h(F^{n-1})=V^{n-1}. 
$$
From this it follows that the sheaves 
${\cal O}_{M}\otimes_{\Co}{\bf F_{n-1}}$ and 
${\cal O}_{M}\otimes_{\Co}{\bf F^{n-1}}$ determine
$V_{n-1}$ and $V^{n-1}$, respectively (see Section 3.1).
By Proposition \ref{Te12} there exists an element $\delta$ of 
$H^{1}(\!M,{\cal O}_{M}\otimes_{\Co}Hom(\!{\bf F_{n-1},F^{n-1}}\!))$ which
determines $V_{n}$. Since the flat bundle 
$Hom(F_{n-1},F^{n-1})$ is isomorphic to $(F_{n-1})^{*}\otimes F^{n-1}$ and 
therefore belongs to ${\cal U}$, we can apply Lemma \ref{Te21}. By the  
lemma there exists an element 
$\beta\in {\em [}\overline{\Omega}_{d}^{1}(Hom(F_{n-1},F^{n-1})){\em ]}
\subseteq H^{1}(M,Hom({\bf F_{n-1},F^{n-1}}))$ such that 
$$
\Pi_{0,1}(\beta)=\delta \ \ \ and \ \ \ \Pi_{1,0}(\beta)=0.
$$
Moreover, $\beta$ defines an extension $F_{n}$ of $F_{n-1}$ by $F^{n-1}$ by 
Proposition \ref{Te12}. From these two statements and Proposition \ref{Te15} 
we conclude that 
$$
h(F_{n})=V_{n} \ \ \ and \ \ \ 
\overline{h}(F_{n})=\overline{h}(F_{n-1})\bigoplus\overline{h}(F^{n-1}).
$$
But $F_{n-1}\in {\cal U}$ by the induction hypothesis and therefore the 
latter direct sum equals 
$$
\bigoplus_{k=1}^{n-1}\overline{h}(F^{k})\bigoplus 
\overline{h}(F_{1})=Gr^{*}\overline{h}(F_{n}).
$$
Thus $F_{n}$ belongs to ${\cal U}$. Moreover, $F_{n}$ is an extension of 
the bundle $F_{n-1}$ by the bundle $F^{n-1}$ and by the induction hypothesis 
these two bundles
are $T_{n-1}^{u}$- and $T_{1}^{u}$-topologically trivial, respectively. So, 
$F_{n}$ is $T_{n}^{u}$-topologically trivial.\\ 
The proof is complete.\ \ \ \ \ $\Box$  
\begin{R}\label{nil}
{\rm If in Theorem \ref{Te1} the form $\omega_{2}$ is nilpotent then it is
$\overline{\partial}$-gauge equivalent to an antiholomorphic nilpotent form.}
\end{R}
\sect{\hspace{-1em}. Proof of Theorem \ref{maintheo}.}
To prove the theorem we have to prove propositions of Section 2.
\\
{\bf Proof of Proposition \ref{Prop1}.}(1)\ 
Let $\psi$ be a diagonal $\overline{\partial}$-closed (0,1)-form on $M$.
According to the Hodge decomposition
\begin{equation}\label{E}
\psi=\widetilde{\psi}+\overline{\partial}f,
\end{equation}
where $\widetilde{\psi}$ is a diagonal, harmonic (0,1)-form. 
Put $h_{\psi}:=\exp(f).$
Then we have $d_{h_{\psi}}(\omega)\in {\cal E}_{\widetilde{\psi}}$
for any $\omega\in {\cal E}_{\psi}$.\ \ \ \ \ $\Box$\\
(2)\
Let $\psi$ be a diagonal, harmonic (0,1)-form on a compact K\"{a}hler 
manifold $M$ (which, in particular, is $d$-closed antiholomorphic). 
Then we determine a flat vector bundle $E_{\psi}$ over $M$
as $=i_{d}(\psi-\overline{\psi})$. As it follows from arguments used in the 
proof of Theorem \ref{Te1},
$E_{\psi}\in {\cal U}_{\oplus}^{n}$, i.e., it is a direct sum of 
rank 1 topologically trivial flat vector bundles with structure group
$U_{1}(\Co)$. Moreover, each element of ${\cal U}_{\oplus}^{n}$ coincides
with $i_{d}(\psi-\overline{\psi})$ for some diagonal harmonic (0,1)-form
$\psi$.
This proves part $(a)$.\\
Let now $\omega\in {\cal E}_{\psi}$, i.e. it has a triangular (0,1)-component
$\omega_{2}$ such that $diag(\omega_{2})=\psi$ and satisfies (\ref{form2}).
Furhter, define the mapping $\tau_{\psi}$ by
\begin{equation}\label{tau}
\tau_{\psi}(\omega ):=\omega-(\psi-\overline{\psi}).
\end{equation}
The latter form can be thought of as a 1-form with values in 
the flat vector bundle $End(E_{\psi})$ whose
(0,1)-component is nilpotent. In fact,
let $\{g_{i}\}_{i\in I}$ be a family of invertible diagonal matrix functions 
defined on an open covering $\{U_{i}\}_{i\in I}$ and satisfying
$$
dg_{i}=(\psi-\overline{\psi})g_{i},\ \ \ \ \ (i\in I).
$$
Then in a flat coordinate system on $End(E_{\psi})$ form
$\tau_{\psi}(\omega)$ is given by the family 
$\{\theta_{i}:=g_{i}^{-1}\tau_{\psi}(\omega)g_{i}\}_{i\in I}$. Clearly
(0,1)-component of $\theta_{i}$ is nilpotent and therefore
$\tau_{\psi}(\omega)$ is, by definition, $End(E_{\psi})$-valued form with a 
nilpotent (0,1)-component. Simple calculation based on the identities 
$$
d\omega-\omega\wedge\omega=0\ \ {\em and}\ \
d(\psi-\overline{\psi})=
(\psi-\overline{\psi})\wedge (\psi-\overline{\psi})=0
$$
and diagonality of $g_{i}$ and $\psi$ get
$$
d\theta_{i}-\theta_{i}\wedge\theta_{i}=0,\ \ \ \ \ (i\in I).
$$
This proves part $(c)$.\\
Let $h\in C^{\infty}(M,T_{n}(\Co))$. Then it determines an element
from the group $Aut_{\infty}^{t}(E_{\psi})$ of triangular authomorphisms of 
$E_{\psi}$ given by the family $\{h_{i}:=g_{i}^{-1}hg_{i}\}_{i\in I}$.
Substituting these expressions in the definition of the $d$-gauge transform
$d_{h}^{E_{\psi}}$ and taking into account diagonality of $g_{i}$ and
$\psi$ we obtain  $\tau_{\psi}\circ d_{h}=d_{h}^{E_{\psi}}\circ\tau_{\psi}$.
This proves part $(b)$.\\
To finish the proof of proposition observe that the mapping $\tau_{\psi}$
defined on ${\cal E}_{\psi}$ by (\ref{tau}) is injective and has the inverse
defined on the set of $End(E_{\psi})$-valued 1-forms with nilpotent 
(0,1)-components satisfying (\ref{form2}).
\ \ \ \ \ $\Box$
\\
{\bf Proof of Proposition \ref{Tebun}.}
In order to prove the proposition we make use of
the relation between elements of ${\cal E}_{\psi}$  
with a diagonal harmonic (0,1)-form $\psi$ and
$End(E_{\psi})$-valued locally solvable equations
with nilpotent (0,1)-components (see Proposition \ref{Prop1}). 

Let $\omega\in {\cal E}_{\psi}$ and $\eta:=\tau_{\psi}(\omega)$ be
an $End(E_{\psi})$-valued differential
1-form with a nilpotent (0,1)-component satisfying the analog of
(\ref{form2}). As follows from Theorem \ref{Te1}, $\omega$ 
can be reduced by a $d$-gauge transform $d_{g}$ with 
$g\in C^{\infty}(M,T_{n}(\Co))$ to
a form $\omega'\in {\cal E}_{\psi}$ with the type decomposition 
$\omega_{1}'+\omega_{2}'$ such that
$\omega_{2}'-diag(\overline{\omega_{2}'})\in {\cal E}_{\psi}$.
Set now 
$$
\widetilde{\eta}_{1}:=\tau_{\psi}(\omega_{1}'+diag(\omega_{2}'))\ \
and\ \
\widetilde{\eta}_{2}:=\tau_{\psi}(\omega_{2}'-
diag(\overline{\omega_{2}'})).
$$
Then clearly $\tau_{\psi}(\omega_{1}'+\omega_{2}')=
\widetilde{\eta}_{1}+\widetilde{\eta}_{2}$ (type decomposition).
According to Proposition \ref{Prop1} $(2b)$, 
$d_{g}^{E_{\psi}}(\eta)=\tau_{\psi}(\omega')
=\widetilde{\eta}_{1}+\widetilde{\eta}_{2}$, where 
$g$ is thought of now as an element of $Aut_{\infty}^{t}(E_{\psi})$.
It remains to prove that $\partial{\widetilde{\eta}_{1}}=0$ and 
$\widetilde{\eta}_{2}$ is a $d$-closed antiholomorphic 1-form.

Actually, (0,1)-form $\widetilde{\eta}_{2}$ satisfies, 
by definition,
$End(E_{\psi})$-valued equation (\ref{form2}). This implies 
immediately that it is antiholomorphic and, due to the Hodge 
decomposition (see Section 3.4), is $d$-closed. 
Prove now that $\widetilde{\eta}_{1}$ is $\partial$-closed.
To accomplish this we observe that conditions (\ref{form2}) for forms 
$\omega_{1}'+\omega_{2}'$ and $\omega_{2}'-diag(\overline{\omega_{2}'})$
include, in particular, the following identities
\begin{equation}\label{om1}
\overline{\partial}\omega_{1}'=\omega_{1}'\wedge\omega_{2}'+
\omega_{2}'\wedge\omega_{1}'-\partial\omega_{2}';
\end{equation}
\begin{equation}\label{om2}
\overline{\partial}(-diag(\overline{\omega_{2}'}))
=(-diag(\overline{\omega_{2}'}))\wedge\omega_{2}'+
\omega_{2}'\wedge (-diag(\overline{\omega_{2}'}))-
\partial\omega_{2}'.
\end{equation}
Then subtracting from the first equation the second one we obtain
\begin{equation}\label{om1hol}
\overline{\partial}(\omega_{1}'+diag(\overline{\omega_{2}'}))=
(\omega_{1}'+diag(\overline{\omega_{2}'}))\wedge\omega_{2}'+
\omega_{2}'\wedge (\omega_{1}'+diag(\overline{\omega_{2}'})).
\end{equation}
Let us consider now the flat vector bundle 
$F:=i_{d}(\omega_{2}'-diag(\overline{\omega_{2}'}))$. Then 
from equation (\ref{om1hol}) it follows that 
$\omega_{1}'+diag(\overline{\omega_{2}'})$ is an $End(F)$-valued holomorphic
1-form (see Section 3.5). Since $F$ belongs to the class ${\cal U}$ which is
closed with respect to tensor products and duality, Lemma
\ref{Te21} (b) implies in this case that 
$\omega_{1}'+diag(\overline{\omega_{2}'})$
is a $d$-closed $End(F)$-valued form. But by the definition $End(F)$
is antiholomorphically isomorphic to $End(Gr^{*}F)$ which is,
in turn, coincides with $End (E_{\psi})$
(see the proof of Proposition \ref{Prop1}). This 
shows that $\omega_{1}'+diag(\overline{\omega_{2}'})$ regarding now as
an $End(E_{\psi})$-valued 1-form is $\partial$-closed. It remains to note that
the latter form coincides with  $\widetilde{\eta}_{1}$. The proof 
of Proposition \ref{Tebun} is complete.\ \ \ $\Box$

Let now $\{\widetilde{\eta}_{1}\}$ and $\{\widetilde{\eta}_{2}\}$ be the harmonic components in the 
Hodge decomposition of $End(E_{\psi})$-valued forms 
$\widetilde{\eta}_{1}$ and $\widetilde{\eta}_{2}$,
respectively. Then $End(E_{\psi})$-valued condition (\ref{form2}) and 
$\partial\overline{\partial}$-lemma of Section 3.4 get
$$
[\{\widetilde{\eta}_{i}\},\{\widetilde{\eta}_{i}\}]=0,\ i=1,2;
$$
$$
[\{\widetilde{\eta}_{1}\},\{\widetilde{\eta}_{2}\}]\ \ {\rm represents}\ \
0\ \ {\rm in}\ \ H^{2}(M,{\bf End(E_{\psi})}).
$$
{\bf Proof of Proposition \ref{Te4}}. 
Let $\alpha_{1},\beta_{1}$ be $End(E_{\psi})$-valued (1,0)-forms, and let
$\alpha_{2},\beta_{2}$ be $End(E_{\psi})$-valued $d$-closed nilpotent 
(0,1)-forms. 
Recall that $E_{\psi}$ is a direct sum of rank 1 topologically trivial flat 
vector bundles with unitary structure group. 
Suppose that
\begin{equation}\label{equi}
d_{g}^{E_{\psi}}(\alpha_{1}+\alpha_{2})=\beta_{1}+\beta_{2}
\end{equation}
for some $C^{\infty}$-authomorphism $g$ of $E_{\psi}$ and  
$\alpha_{1}+\alpha_{2}$ and $\beta_{1}+\beta_{2}$ 
belong to $\tau_{\psi}({\cal E}_{\psi})$.
We have to prove that $g$ is flat.
According to Propositions \ref{Prop1} and \ref{Tebun} there
exist triangular (0,1)-forms $\theta_{1}$ and $\theta_{2}$ such that\\
$(i)\ \ \tau_{\psi}(\theta_{1}-\overline{\psi})=\alpha_{2},\ \ 
\tau_{\psi}(\theta_{2}-\overline{\psi})=\beta_{2};$\\
$(ii)\ \ \ diag(\theta_{i})=\psi,\  \ (i=1,2);$\\
$(iii)\ \ \theta_{i}-diag(\overline{\theta_{i}})\in {\cal E}_{\psi}$,
\ $(i=1,2)$.\\
If we now identify the group of $C^{\infty}$-authomorphisms of $E_{\psi}$
with $C^{\infty}(M,GL_{n}(\Co))$ (as in the case of triangular authomorphisms)
then arguing as in the proof of Proposition \ref{Prop1} we obtain
$\tau_{\psi}\circ d_{g}=d_{g}^{E_{\psi}}\circ\tau_{\psi}$. In particular,
(\ref{equi}) implies
$$
\overline{\partial}g=\theta_{1}g-\theta_{2}g.
$$
But this is a special case of equation (\ref{form15}). Applying 
Proposition \ref{Te20} we conclude that $g$ is a holomorphic  section of 
flat vector bundle 
$V:=Hom(i_{d}(\theta_{2}-diag(\overline{\theta_{2}})),\penalty-10000
i_{d}(\theta_{1}-diag(\overline{\theta_{1}})))$. 
This vector bundle belongs to the class ${\cal U}$, and therefore $g$ is 
$d$-closed by Lemma \ref{Te21}(b). Since by definition $V$ is
antiholomorphically isomorphic to $End(E_{\psi})$, the authomorphism 
$g$ of $E_{\psi}$ is $\partial$-closed.
Applying now the Hodge decomposition of Section 3.4 we deduce that $g$ is 
locally constant, i.e., flat.
\ \ \ \ \ $\Box$ 
\\
{\bf Proof of Proposition \ref{Te3}}.  Let 
$\alpha$ be a holomorphic 
$End(E_{\psi})$-valued form and $\theta$ be an 
antiholomorphic nilpotent one and
let the 2-form 
$[\alpha+\theta,\alpha+\theta]$ represent $0$ in
$H^{2}(M,{\bf End(E_{\psi})})$. 
We have to prove that there exists a 
section $h$, unique up to an additive flat summand, such that the equation 
$$
df=(\alpha+\theta+\partial h)f
$$
is locally solvable.
To accomplish this we first remark that 
$\partial\overline{\partial}$-lemma of
Section 3.4 implies that
$$
\alpha\wedge\theta+\theta\wedge\alpha=[\alpha+\theta,\alpha+\theta]
$$
since the form on the right represents 0 in $H^{2}(M,{\bf End(E_{\psi})})$.
Applying the $\partial\overline{\partial}$-lemma to the left-hand side
and taking into account the holomorphicity
of $\alpha$ we then obtain
\begin{equation}\label{form21}
\overline{\partial}\alpha-
\alpha\wedge\theta-\theta\wedge\alpha=\partial\overline{\partial}P
\end{equation}
for some $C^{\infty}$-section $P$ of $End(E_{\psi})$. 
Since according to our assumption $d\theta-\theta\wedge\theta=0$, the 
arguments similar to those of
Proposition \ref{Prop1} show that there exists a triangular (0,1)-form
$\eta$ defined on $M$ such that 
$\eta-diag(\overline{\eta})$ satisfies (\ref{form2}),
$diag(\eta)=\psi$ and 
$\tau_{\psi}(\eta-diag(\overline{\eta}))=\theta$.
Then in the global $C^{\infty}$-coordinates on $End(E_{\psi})$ (chosen as
in the proof of Proposition \ref{Prop1}) (\ref{form21}) 
can be written as
$$
\overline{\partial}\alpha'-\eta\wedge\alpha'-\alpha'\wedge\eta=
\partial\beta+
diag(\overline{\eta})\wedge\beta+
\beta\wedge diag(\overline{\eta}),
$$
(see Section 3.5). Here $\alpha:=g_{i}^{-1}\alpha'g_{i}$ on $U_{i}$ and
$\overline{\partial}P:=g_{i}^{-1}\beta g_{i}$ on $U_{i}$ and 
$\{g_{i}\}_{i\in I}$ is a family of
invertible diagonal matrix function satisfying
$dg_{i}=(\psi-\overline{\psi})g_{i}$ on $U_{i}$.
Consider now flat vector bundle $F:=i_{d}(\eta-diag(\overline{\eta}))$ of the
class ${\cal U}$. If we now think of $\alpha$ as
an $End(F)$-valued (1,0)-form ($F$ is $C^{\infty}$-trivial)
then the left-hand side of the previous expression determines its
$\overline{\partial}$-differential.
But the right-hand side  shows that 
$\overline{\partial}\alpha$ is a $\partial$-exact
$End(F)$-valued form.  
The proof of the theorem will be complete if we find a 
$C^{\infty}$-section 
$h$ of $End(F)$ such that $\alpha+\partial h$ is a holomorphic 
$End(F)$-valued form. 
Actually, let $\{f_{i}\}_{i\in I}$ be a family of triangular invertible
$C^{\infty}$-functions determined on the open covering $\{U_{i}\}_{i\in I}$ 
(the same covering as for $\{g_{i}\}_{i\in I}$ above) and satisfying
$df_{i}=(\eta-diag(\overline{\eta}))f_{i}$, $(i\in I)$.
Then the holomorphicity of $\alpha+\partial h$ is equivalent to the 
equation 
$$
\overline{\partial}(\alpha'+\gamma)=(\alpha'+\gamma)\wedge\eta+
\eta\wedge(\alpha'+\gamma),
$$
where $\gamma=f_{i}\partial hf_{i}^{-1}$ on $U_{i}$.
The latter  equation, in turn, determines the $End(E_{\psi})$-valued equation
\begin{equation}\label{form22}
\overline{\partial}(\alpha+\widetilde\gamma)=
(\alpha+\widetilde\gamma)\wedge\theta+\theta\wedge (\alpha+\widetilde\gamma).
\end{equation}
Here $\widetilde\gamma=g_{i}^{-1}\gamma g_{i}$ on $U_{i}$.
Clearly, $\partial (g_{i}^{-1}f_{i})=0$ and therefore
$\widetilde\gamma=\partial (g_{i}^{-1}f_{i}hf_{i}^{-1}g_{i})$ on $U_{i}$.
But $\{g_{i}^{-1}f_{i}hf_{i}^{-1}g_{i}\}_{i\in I}$ determines a section
$\widetilde h$ of $End(E_{\psi})$. So $\widetilde\gamma=\partial\widetilde h$.
Equation (\ref{form22}) is one of the conditions of local solvability
containing in (\ref{form2}).
One observes that (\ref{form2}) in our case is equivalent to 
the fulfillment of (\ref{form22}) together with the identity 
\begin{equation}\label{hol}
(\alpha+\partial\widetilde h)\wedge (\alpha+\partial\widetilde h)=0,
\end{equation}
since $\partial(\alpha+\partial\widetilde h)=0$ by 
assumptions of the proposition.
To check this identity we first note that $End(E_{\psi})$
is antiholomorphically isomorphic to $End(F)$. This isomorphism is given 
locally by conjugations by matrix functions $f_{i}^{-1}g_{i}$ \ $(i\in I)$ and 
so it commutes with the operator $\wedge$. 
Therefore it suffices to prove an identity similar to (\ref{hol}) for 
$\alpha+\partial h$. Here $\alpha$ is thought of as an 
$End(F)$-valued section (image of $\alpha$ by the above isomorphism). 
Furthermore, because $\alpha\wedge\alpha=0$ one has
\begin{equation}\label{form23}
(\alpha+\partial h)\wedge(\alpha+\partial h)=
\partial(h\alpha-\alpha h+h\partial h).
\end{equation}
This implies that $End(F)$-valued holomorphic 1-form $\alpha+\partial h$
is $\partial$-exact. 
Applying now Lemma \ref{Te21}(b) to this form one concludes that the identity
(\ref{hol}) holds. The uniqueness part of the proposition follows from the
fact that there is a unique up to a flat additive summand section $h$
such that $\alpha+\partial h$ is $End(F)$-valued holomorphic (see Lemma
\ref{Te21}(b)).

Thus it remains to find the section $h$ such that $\alpha+\partial h$ is a
holomorphic $End(F)$-valued 1-form.
We do this by a procedure reducing the $n$-dimensional statement to 
the $(n-1)$-dimensional one; here $n$ is the dimension of a fibre of $End(F)$. \\
We begin with the following remark.
Since $End(F)\in {\cal U}$ it can be regarded as an extension of 
a rank 1 flat vector bundle $F_{1}$ with unitary structure group
by a flat vector bundle $F_{n-1}\in {\cal U}$, i.e., the following sequence of 
flat vector bundles
$$
0\longrightarrow F_{n-1}\stackrel{i}{\longrightarrow}End(F) 
\stackrel{j}{\longrightarrow}F_{1}\longrightarrow 0
$$
is exact. We can analogously represent  $F_{n-1}$ as an extension of 
a rank 1 flat vector bundle with unitary structure group by a flat 
vector bundle $F_{n-2}\in {\cal U}$ and so on. In particular, $F_{0}$ is a
vector bundle over $M$ with null-dimensional fibre.\\
In the next part of the proof we let the same letters $i,j$ 
denote the corresponding
mappings induced by $i,j$ on the space of differential forms.\\
Let us consider now the $F_{1}$-valued $\partial$-exact (1,1)-form 
$j(\overline{\partial}\alpha)=\overline{\partial}(j(\alpha))$. 
Since $\partial\alpha=0$, we get $\overline{\partial}\alpha=d\alpha$, and hence
$j(\overline{\partial}\alpha)$ is a $d$-exact $F_{1}$-valued 1-form. 
The $\partial\overline{\partial}$-lemma implies then that 
$$
\overline{\partial}(j(\alpha))=\overline{\partial}\partial(g) 
$$
for $g\in C^{\infty}(F_{1})$.
Since $End(F)$ is a trivial extension of $F_{1}$ by $F_{n-1}$ in the category 
of $C^{\infty}$-bundles, there exists an $End(F)$-valued 
$C^{\infty}$-section $k_{1}$ 
such that $j(k_{1})=g$. If we put now $\alpha_{1}:=\alpha-\partial k_{1}$ then 
$$
\partial\alpha_{1}=\partial\alpha=0 \ \ \ and \ \ \ 
\overline{\partial}(j(\alpha_{1}))=j(\overline{\partial}\alpha_{1})=
j(\overline{\partial}\alpha-\overline{\partial}\partial k_{1})=0.
$$
It follows from the second identity  that $\overline{\partial}\alpha_{1}$
can be regarded as an $F_{n-1}$-valued form.
Since $End(F)=F_{1}\oplus F_{n-1}$ in the class of antiholomorphic vector 
bundles, the mapping 
$$
i:H^{1}(M,\overline{\Omega}^{1}(F_{n-1}))\longrightarrow 
H^{1}(M,\overline{\Omega}^{1}(End(F)))
$$ 
is an injection. 
Moreover,  the $\partial$-exactness of $\overline{\partial}\alpha$ implies that
$$
i([\overline{\partial}\alpha_{1}])=0\in 
H^{1}(M,\overline{\Omega}^{1}(End(F)))
$$ 
and therefore 
$[\overline{\partial}\alpha_{1}]=0\in 
H^{1}(M,\overline{\Omega}^{1}(F_{n-1}))$.
From this it follows that $\overline{\partial}\alpha_{1}$ is an 
$F_{n-1}$-valued $\partial$-exact form. 
Starting with the $F_{n-1}$-valued form $\overline{\partial}\alpha_{1}$ and 
proceeding in the same way we can now find a 
$C^{\infty}$-section $k_{2}$ such that 
for
$$
\alpha_{2}:=\alpha_{1}-\partial k_{2}
=\alpha-\partial k_{1}-\partial k_{2}
$$ 
$\overline{\partial}\alpha_{2}$ is an $F_{n-2}$-valued $\partial$-exact 
(1,1)-form.\\ 
Continuing in this fashion we obtain after $n$ steps the form
$\alpha_{n}:=\alpha_{n-1}-\partial k_{n-1}$ such that 
$\overline{\partial}\alpha_{n}$ is an $F_{0}$-valued 
$\partial$-exact (1,0)-form, i.e., $\overline{\partial}\alpha_{n}=0$. 
If we now set
$$
h:=-\sum_{i=1}^{n}k_{i}
$$ 
then $\alpha+\partial h$ equals 
the holomorphic $End(F)$-valued 1-form $\alpha_{n}$. \ \ \ \ \ $\Box$
\\ 
\begin{R}\label{Te23}
{\rm 
If the form $\alpha$ of Proposition \ref{Te3} is, in addition, 
triangular then the section $h$ can also be chosen as triangular. 

In fact, let $t_{n}$ be the Lie algebra of the 
Lie group $T_{n}(\Co)$ of upper triangular matrices. 
The vector space $t_{n}$ is invariant with respect to the linear operators 
$(A^{t})^{-1}\otimes A:M_{n}(\Co)\longrightarrow M_{n}(\Co)$ 
with $A\in T_{n}(\Co)$. Therefore  
there exists a sub-bundle $T\in {\cal U}$ of the bundle
$End(F)$ of the proof of Proposition \ref{Te3} with a fibre isomorphic to 
$t_{n}$. In fact, 
the latter bundle is defined by a cocycle of the form
$$
\{(c_{ij}^{t})^{-1}\otimes c_{ij};\ \ c_{ij}\in T_{n}(\Co)\}
$$ 
(see Section 3.5).
Since the form $\alpha$ is $t_{n}$-valued and $T\in {\cal U}$, this form is 
a $T$-valued (1,0)- one. We can now apply the arguments of the proof
of Theorem \ref{Te3} to $\alpha$ but with $T$ instead of
$End(F)$.
In this way we obtain the required section $h$ but in this case with values 
in $T_{n}(\Co)$.}
\end{R}
\sect{\hspace{-1em}. Proof of Theorems \ref{te5}, \ref{te6}.}
{\bf Proof of Theorem \ref{te5}.}
Let $V_{2}(M)$ be a class of flat vector bundles over $M$ whose elements
are constructed by homomorphisms from $S_{2}^{u}(M)$.
According to the assumptions, for any $E\in V_{2}(M_{1})$ 
there exists $F\in V_{2}(M_{2})$ such that $f^{*}F\cong E$.
Moreover every such bundle $E$ is, by definition, determined by 
$Gr^{*}E=E_{1}\oplus E_{2}$ 
and an element of $H^{1}(M_{1}, {\bf Hom(E_{2},E_{1})})$. Here $E_{1},E_{2}$
are topologically trivial rank 1 flat vector
bundles with unitary structure group. Then the 
conditions of the theorem imply\\
{\bf Statement.}
{\em
For every topologically trivial rank 1 flat vector bundle $V_{1}$ over $M_{1}$
with unitary structure group there exists a topologically trivial flat vector 
bundle $V_{2}$ over $M_{2}$ with unitary structure group such that 
$f^{*}V_{2}=V_{1}$ and
$f^{*}(H^{1}(M_{2}, {\bf V_{2}}))=H^{1}(M_{1}, {\bf V_{1}})$. }

Let now $\rho:\pi_{1}(M_{1})\longrightarrow GL_{n}(\Co)$ be a homomorphism 
of the class $S_{n}(M_{1})$. Then according to Theorem \ref{maintheo},
$\rho$ is uniquely defined by $End(E')$-valued harmonic (1,0)-form 
$\alpha$ and harmonic nilpotent (0,1)-form $\eta$ satisfying 
$[\alpha+\eta,\alpha+\eta]$ represents 0 in
$H^{2}(M_{1}, {\bf End(E')})$. Here $E'$ is a direct sum of
topologically trivial 
rank 1 flat vector bundles with unitary structure group. Furthemore, from the 
above Statement it follows that there exist a flat vector bundle $F'$ over
$M_{2}$ isomorphic to a direct sum of topologically trivial rank 1 flat 
vector bundles with unitary
structure group and $End(F')$-valued harmonic (1,0)-form $\alpha'$ and
harmonic nilpotent (0,1)-form $\eta'$ such that 
$$
f^{*}End(F')=End(E'),\ \ f^{*}(\alpha')=\alpha,\ \ f^{*}(\eta')=\eta.
$$
In addition, assume that $[\alpha',\eta']$ represents 0 in
$H^{2}(M_{2}, {\bf End(F')})$. The above conditions imply also
$[\alpha',\alpha']=[\eta',\eta']=0$ and therefore
the triple $(End(F'),\alpha',\eta')$ determines a 
representation $\rho'\in S_{n}(M_{2})$. Then the uniqueness part of 
Theorem \ref{maintheo} (see Proposition \ref{Te3}) yields 
$\rho=\rho'\circ f_{*}$.

Thus it remains to prove that $[\alpha',\eta']$ represents 0 in
$H^{2}(M_{2}, {\bf End(F')})$. Note that $f^{*}([\alpha',\eta'])=
[\alpha,\eta]$ represents 0 in $H^{2}(M_{1}, {\bf End(E')})$.
The required statement then is a consequence of the following general result.
 
Let $N\stackrel{f}{\longrightarrow}M$ be a surjective mapping of compact 
K\"{a}hler manifolds, $E$ be a flat vector bundle over
$M$ with unitary structure group. 
\begin{Proposition}\label{null}
Let $\alpha\in {\cal E}^{1,1}(E)$ be a $d$-closed $E$-valued  
form. If $f^{*}(\alpha)\in {\cal E}^{1,1}(f^{*}E)$ is $d$-exact then 
$\alpha$ is also $d$-exact.
\end{Proposition}
{\bf Proof.}
Consider the flat vector bundle $f^{*}E$ over $N$ and 
$d$-exact form $f^{*}(\alpha)\in {\cal E}^{1,1}(f^{*}E)$.
Since this bundle has unitary structure group and $N$ is a compact K\"{a}hler
manifold, there exists 
$h\in C^{\infty}(f^{*}E)$ such that 
$f^{*}(\alpha)=\overline{\partial}\partial h$. 
Let $N\stackrel{p_{1}}{\longrightarrow}Y\stackrel{p_{2}}{\longrightarrow}M$
be the Stein factorization of $f$. Here the fibres of $p_{1}$ are connected
and $p_{2}$ is a finite analytic covering. For a point
$x\in M$ consider an open neighborhood $U_{x}$ of $x$ such that 
$E\mid_{U_{x}}$ is the trivial flat vector bundle. Then $f^{*}E$ is trivial 
over $f^{-1}(U_{x})$ and for any fibre $V$ of $f$ 
over a point of $U_{x}$ the restriction 
$\alpha_{V}:=f^{*}(\alpha)|_{V}=0$. This implies that
$h|_{V}$ is locally constant. (To prove this fact in the case of
singular $V$ one has to pull back $h$ 
to its desingularisation.) 
Then there exists a section $h'$ of
$p_{2}^{*}E$ such that $p_{2}^{*}\alpha=\overline{\partial}\partial h'$ 
on non-singular part of $Y$. Consider now the average of $h'$ 
over points of regular fibres of $p_{2}$
$$
h''(y):=\frac{1}{\#\{p_{2}^{-1}(y)\}}\sum_{z\in p_{2}^{-1}(y)}h'(z),\ \ \  
(y\in M).
$$
Clearly $h''$ is a bounded section of $E$ smooth at regular values
of $p_{2}$. 
Then  $\alpha=\overline{\partial}\partial h''$ outside of a proper 
analytic subset of $M$. Moreover, according to assumptions of
the proposition $\alpha$ is locally  
$\overline{\partial}\partial$-exact. Further, boundedness of $h''$
together with regularity of the operator $\partial\overline{\partial}$ imply 
that $h''$ can be extended to $M$ as a $C^{\infty}$-section of $E$
satisfying $\alpha=\overline{\partial}\partial h''$. This shows
that $\alpha$ is $d$-exact.\ \ \ $\Box$

The proof of Theorem \ref{te5} is complete.\ \ \ $\Box$.\\
\\
{\bf Proof of Theorem \ref{te6}.}
Let $\tau:\pi_{1}(M_{1})\longrightarrow T_{2}^{u}$ be a representation
of the class $S_{2}^{u}(M_{1})$.
Clearly $Ker(\tau)$ contains $\pi_{1}(M_{1})''$ and so $\tau$ determines a
homomorphism $\penalty -10000$
$\tau_{1}:\pi_{1}(M_{1})/\pi_{1}(M_{1})''\longrightarrow 
T_{2}^{u}$. Furthemore, according to the assumptions of the theorem, there 
exists a homomorphism 
$\tau_{2}:\pi_{1}(M_{2})/\pi_{1}(M_{2})''\longrightarrow
T_{2}^{u}$ such that $\tau_{1}=\tau_{2}\circ f_{*}$ whose diagonal elements
have the logarithm. Obviously, we can 
extend $\tau_{2}$ to a homomorphism $\tau':\pi_{1}(M_{2})\longrightarrow
T_{2}^{u}$ of the class $S_{2}^{u}(M_{2})$ satisfying $\tau=\tau'\circ f_{*}$. 
Thus the 
conditions of Theorem \ref{te5} are fulfilled. According to this theorem 
for any representation $\rho:\pi_{1}(M_{1})\longrightarrow GL_{n}(\Co)$ of
the class $S_{n}(M_{1})$ there exists a representation
$\rho':\pi_{1}(M_{2})\longrightarrow GL_{n}(\Co)$ such that
$\rho=\rho'\circ f_{*}$. The latter, in particular, shows that 
$Ker f_{*}$ belongs to the kernel of every matrix representation of the class 
$S_{n}(M_{1})$ $(n\geq 1)$. But by the assumption of the theorem 
$\pi_{1}(M_{1})$ belongs to the class $S$. Therefore $Ker f_{*}=\{e\}$ and
$f_{*}$ is an injective homomorphism. Moreover, from the Stein factorization
of $f$ one obtains that $f_{*}(\pi_{1}(M_{1}))$ is a subgroup of a
finite index in  $\pi_{1}(M_{2})$.\ \ \ $\Box$
\sect{\hspace{-1em}. Concluding Remarks.}
All results of this paper hold also true for the class of manifolds dominated
by a compact K\"{a}hler. We recall the following
\begin{D}\label{domin}
A manifold $M$ is said to be dominated by a compact K\"{a}hler manifold
$N$ if there exists a complex surjective mapping $f:N\longrightarrow M$.
\end{D}
Let M be a manifold dominated by a compact K\"{a}hler manifold
$N$: $N\stackrel{f}{\longrightarrow}M$ and $E$ be a flat vector bundle over
$M$ with unitary structure group. The proof of the following proposition
is similar to that of Proposition \ref{null}.
\begin{Proposition}\label{formality}
(a)\ \ \ Let $\alpha\in {\cal E}^{0,1}(E)$ be an $E$-valued  
$\overline{\partial}$-closed (0,1)-form. Then there exists a 
$C^{\infty}$-section $h$ of $E$ such that $\alpha-\overline{\partial}h$ is
$d$-closed.\\
(b)\ \ \ Let $E$-valued (1,1)-form $\beta$ satisfy $d\beta=0$  and
$\beta=\partial\gamma$ for some $E$-valued (0,1)-form $\gamma$. Then 
there exists an $E$-valued function $g$ such that 
$\beta=\partial\overline{\partial}g$.
\end{Proposition}
Using this result and applying the very same arguments one can prove 
the validity of the results of the paper for the class of manifolds 
dominated by compact K\"{a}hler ones.        

\noindent Department of Mathematics\\
University of Toronto, Toronto, Ontario M5S 1A1 CANADA
\end{document}